%% file: Robust.tex
\documentclass[11pt]{article}
\usepackage{amsmath,amsthm,amsfonts,amssymb,amscd}
\usepackage{enumitem}
\usepackage{color}
\usepackage{float}
\usepackage{soul}
\usepackage{graphicx}
\usepackage{scalerel}
\usepackage{tikz}
\usepackage{mathpazo}
\usepackage{fancyhdr}
\usepackage{empheq}
\usepackage[margin=1in]{geometry}
\usepackage{mysty}
\usepackage[colorlinks=true]{hyperref}
\hypersetup{  
   bookmarks=true,
   backref=true,
   pagebackref=false,
   colorlinks=true,
   linkcolor=blue,
   citecolor=red,
   urlcolor=blue,
}
\graphicspath{{images/}}

\begin{document}
\title{\bf Sharp, strong and unique minimizers for low complexity robust recovery}
%\title{\bf Quantitative characterizations of sharp and strong solutions with applications to robust recovery}
\date{}
\author{Jalal Fadili\footnote{Normandie Universit\'e, ENSICAEN, UNICAEN, CNRS, GREYC, France; email:Jalal.Fadili@greyc.ensicaen.fr} \and Tran T. A. Nghia\footnote{Department of Mathematics and Statistics, Oakland University, Rochester, MI 48309, USA; email: nttran@oakland.edu. Research of this author was supported by the US National Science Foundation under grant DMS-1816386.} \and  Trinh T. T. Tran\footnote{Department of Mathematics and Statistics, Oakland University, Rochester, MI 48309, USA; email: thitutrinhtran@oakland.edu}}

\maketitle
\begin{abstract}
{\small   In this paper, we show the important roles of sharp minima and strong minima for robust recovery. We also obtain several characterizations of sharp minima for convex regularized optimization problems. Our characterizations are quantitative and verifiable especially for the case of decomposable norm regularized problems including sparsity, group-sparsity, and low-rank convex problems. For  group-sparsity optimization problems, we show that a unique solution is a strong solution and obtain quantitative characterizations for solution uniqueness.} 
\end{abstract}

\tableofcontents

\input{tex/sec_intro}
\input{tex/sec_prelim}

\input{tex/sec_unique}
\input{tex/sec_sharp}
\input{tex/sec_strong_group}

\input{tex/sec_numerics}
\input{tex/sec_conclusion}

\bibliographystyle{plain}
\bibliography{biblio}

\end{document}

%% file: tex/sec_intro.tex
%%%%%%%%%%%%%%%%%%%%%%%%%%%%%%%%%%%%%%%%%%%%%%%%%%%%%%%%
\section{Introduction}
\label{sec:intro}
%%%%%%%%%%%%%%%%%%%%%%%%%%%%%%%%%%%%%%%%%%%%%%%%%%%%%%%%

%%%%%%%%%%%%%%%%%%%%%%%%%%%%%%%%%%%%%%%%%%%%%%%%%%%%%%%%
\subsection{Problem statement}
Inverse problems and regularization theory is a central theme in various areas of engineering and science. A typical case is where one observes a vector $y_0 \in \R^m$ of linear measurements according to
\begin{equation}\label{y0}
y_0=\Phi x_0 
\end{equation}
where $\Phi \in \R^{m\times n}$ is known, and is typically an idealization of the sensing mechanism in signal/imaging science applications, or the design matrix in a parametric statistical regression problem. The vector $x_0 \in \R^n$ is the unknown quantity of interest. 

Solving the inverse problem associated to the linear forward model \eqref{y0} amounts to recovering $x_0$, either exactly or to a good approximation, knowing $y_0$ and $\Phi$. This is however a quite challenging task especially when the linear system \eqref{y0} is underdetermined with $m \ll n$. In fact, even when $m=n$, $\Phi$ is in general ill-conditioned or even singular. This entails that the linear inverse problem is in general ill-posed.

In order to bring back the problem to the land of well-posedness, it is necessary to restrict the inversion process to a well-chosen subset of $\R^n$ containing the plausible solutions including $x_0$. This can be achieved by adopting a variational framework and solving the following optimization problem
\begin{equation}\label{BP0}
\min_{x\in \R^n}\quad J(x) \qsubjq \Phi x=y_0,
\end{equation}
where $J$ is function which is bounded from below (wlog non-negative). The function $J$, known as a regularizer, is designed in such way that it is the smallest on the sought-after solutions. Popular example in signal/image processing and machine learning is the $\ell_1$ norm to promote sparsity, the $\ell_1/\ell_2$ norm to promote group sparsity, analysis sparsity seminorm (i.e., $J = J_0 \circ D^*$, $D^*: \R^n \to \R^p$ is linear and $J_0$ is the $\ell_1$ or $\ell_1/\ell_2$ norm); see \cite{BDE09,CT04,DH01, JKL15,RRN12, RF08, T96, VPDF13,YL06} to cite a few. Another popular example is the nuclear norm, i.e., the $\ell_1$ norm of the singular values of a matrix, to recover low-rank matrices \cite{CRPW12,CR09}. These prototypical cases are included in the class of decomposable norms \cite{NRWY09,CR13,FPVDS13} that covers a wide range of optimization-based recovery problems in different areas of data science.

When the observation is subject to errors, the system \eqref{y0} is modified to   
\begin{equation}\label{yw}
y=\Phi x_0 +\omega
\end{equation}
where $\omega$  accounts either for noise and/or modeling errors. The errors can be either deterministic (in this case, one typically assumes to know some bound on $\|\omega\|$) or random (in which case its distribution is assumed to be known). Throughout, $\omega$ will be assumed deterministic with $\|\omega\|\le \delta$, with $\delta > 0$ known.

To recover $x_0$ from \eqref{yw}, one solves the noise-aware form (aka Mozorov regularization in the inverse problems literature)
\begin{equation}\label{BP1}
\min_{x \in \R^n} J(x)\quad \mbox{subject to}\quad \|y - \Phi x\| \le \delta
\end{equation}
or its penalized/Lagrangian form (aka Tikhonov regularization in the inverse problems literature)
\begin{equation}\label{Las3}
\min_{x \in \R^n} \frac{1}{2}\|\Phi x-y\|^2+ \mu J(x) .
\end{equation}
with regularization parameter $\mu>0$. It is well-known that \eqref{BP1} and \eqref{Las3} are formally equivalent in the sense that there exists a bijection between $\delta$ and $\mu$ such that the both problems share the same set of solutions.

We say that {\em robust recovery} occurs when solutions to \eqref{BP1} or \eqref{Las3} are close enough to the original signal $x_0$ as long as $\mu$ is appropriately chosen (as a function of the noise level $\delta$). For $\ell_1$-regularized optimization problem, a significant result was achieved in \cite[Theorem~4.7]{GHS11}, which shows that when $\mu$ is proportional to $\delta$, the error between any solution of problem \eqref{Las3} to $x_0$ is $\mathcal{O}(\delta)$\footnote{This justifies the conventional terminology "linear convergence rate".} if and only if $x_0$ is the unique solution to \eqref{BP0}. The structure of the $\ell_1$ norm plays a crucial role in their result. Sufficient (but not necessary in general) conditions for robust recovery with a linear convergence rate have been proposed for general class of regularizers by solving either \eqref{BP1} or \eqref{Las3}; see \cite{CRPW12, FPVDS13, V14, VPF15, Tropp15} and references therein. For instance \cite{CRPW12} and \cite{Tropp15} provide a sufficient condition via the geometric notion of {\em descent cone} for robust recovery with linear rate by solving \eqref{BP1} (but not \eqref{Las3}). Injectivity of $\Phi$ on the descent cone turns out to be a sufficient and necessary condition for solution uniqueness of \eqref{BP0}. However, solution uniqueness is not enough to guarantee robust recovery with linear rate for general regularizers. 

In fact, we observe that \cite{CRPW12} needs $x_0$ to be a {\em sharp} solution to problem \eqref{BP0}. A sharp minimum, introduced by Crome and Polyak \cite{C78,P79,P87}, is a point around which the objective value has a first-order growth (see Section~\ref{sec:prelim} for a rigorous definition and discussion). It is certainly sufficient for solution uniqueness and indeed is called {\em strong uniqueness} in \cite{C78}. In the case of polyhedral problems, e.g., problem \eqref{BP0} with $\ell_1$ or $\ell_\infty$ norms, solution uniqueness is equivalent to sharp minima. This is the key reason why \cite{GHS11} only needs solution uniqueness for robust recovery with linear rate. In the optimization community, sharp minima is well-known for its essence to finite convergence of the proximal point algorithms; see, e.g., \cite{P87, R76}.

Sharp minima is actually hidden in both {\em exact recovery} (by solving \eqref{BP0}) and robust recovery. As we will show later, sharpness is closely related to many other results on exact and robust recovery, for instance those involving conditions based on non-degenerate dual certificates (aka {\em Nondegenerate Source Condition}) and restricted injectivity of $\Phi$, or the {\em Null Space Property} and its variants; see \cite{DH01, CT04, T04, F05, T06, CR09, BDE09, GHS11, FR13, VPDF13, FPVDS13, CR13, VPDF13, V14, VPF15, ZYY16}. When $\Phi$ is drawn from an appropriate random ensemble, existing sample complexiy bounds (i.e., a lower-bound on $m$ depending on the "intrinsic" dimension of the model subspace containing $x_0$ and logarithmically on $n$) rely on verifying the afore-mentioned conditions with high probability on $\Phi$; see \cite{FR13,VPF15,Tropp15} for comprehensive reviews. For Gaussian measurements, sample complexity bounds for exact recovery are provided in terms of the Gaussian width in \cite{CRPW12} and the statistical dimension in \cite{ALMT14}, and the latter has been shown to precise predictions about the quantitative aspects of the phase transition for exact recovery from Gaussian measurements.

Another popular notion in optimization \cite{BS00,RW98} is that of {\em strong minima}, which is a sufficient condition for solution uniqueness. It is also used implicitly in many results on robust recovery. A strong minimum is a point around which the objective value has a quadratic (or second-order) growth. Strong minima are certainly weaker than sharp minima. In the $\ell_1$ case, we observe  that the cost function in \eqref{Las3} belongs to the bigger class of convex {\em piecewise linear-quadratic functions} \cite{RW98}. A unique solution to this function is also a strong solution. This simple observation was used recently in \cite{BLN19} to obtain some characterizations for solution uniqueness to \eqref{Las3} for the case of the Lasso problem. Necessary and sufficient conditions for solution uniqueness to \eqref{Las3} where $J$ is piecewise linear are also studied in \cite{G17,MS19,ZYY16}.
%For $\ell_1$ problems, Fuchs \cite{F05} and Tropp \cite{T04} introduced a synthesis criterion that guarantees solution uniqueness to \eqref{Las3} when $\mu$ is sufficiently small and $y=y_0$. Moreover,  this unique solution has the same sign and support with the original signal $x_0$. It means that one only needs to solve problem \eqref{Las3} with some small  $\mu$ to get enough information about the sign and support of $x_0$. This   criterion has been studied extensively in the literature; see, e.g., \cite{BDE09, FPVDS13, VPDF13, T06}. Solution uniqueness to problem \eqref{Las3} is significant in this direction.  \\

%%%%%%%%%%%%%%%%%%%%%%%%%%%%%%%%%%%%%%%%%%%%%%%%%%%%%%%%
\subsection{Contributions and relation to prior work} 
The chief goal of this paper is to show the intricate and important roles of sharp minima and strong minima in robust recovery. Our main contributions are as follows.

\begin{enumerate}[label=$\bullet$]
\item As discussed above, solution uniqueness and sharp minima play a pivotal role for robust recovery with linear rate. However, except the case where the regularizer $J$ is merely polyhedral, a unique solution may not be sharp. For example, as we will show in our Section~\ref{sec:stronggroup}, a unique solution to the group-sparsity problem \eqref{BP0} is always a strong minimum but not a sharp one. The natural question is therefore whether strong minima could also lead to some meaningful results for robust recovery. Section~\ref{sec:unique} answers this question positively without even needing convexity of $J$, unlike the previous work of \cite{V14,VPF15} which showed that a sharp minimizer implicitly guarantees robust recovery with linear rate, but convexity was essential there. We show that if $x_0$ is a strong solution to \eqref{BP0}, then we have robust recovery by solving \eqref{BP1} or \eqref{Las3} with a convergence rate $\mathcal{O}(\sqrt{\delta})$ (Theorem~\ref{SRob}). This rate can be improved to the linear rate $\mathcal{O}(\delta)$ when $J = \|\cdot\|$. Whether this linear rate for robust recovery can be achieved for other regularizers under the sole strong minimality assumption is an open question. Besides that main result, we also revisit in Theorem~\ref{Rob} robust recovery results when $x_0$ is assumed to be a sharp  minimizer by proving linear convergence and providing more explicit bounds in comparison to the results in \cite{CRPW12, FPVDS13, GHS11, ZYY16}.

\item Our second main contribution, at the heart of Section~\ref{sec:sharp}, is to study necessary and sufficient conditions for the properties of sharp and strong minima. As mentioned above, solution uniqueness is characterized via the descent cone \cite{CRPW12}. But it is hard to compute the descent cone for general functions $J$. On the other hand, sharp and strong minima can be characterized explicitly via first and second-order analysis \cite{BS00,P87,RW98}. We leverage these tools (see Theorem~\ref{mtheo} and Proposition~\ref{Stheo}) to get equivalent charaterizations of unique/sharp/strong solutions, and to obtain a quantitative condition for convex regularized problems that allows us to check numerically that a solution is sharp. This condition bears similarities with the Nondegenerate Source Condition of \cite{GHS11,V14,VPF15}. Our approach does not rely on any polyhedral structure. This distinguishes our results from several ones in the literature as in \cite{G17,MS19,ZYY16}, and allows us to work on more general optimization problems. We also believe that non-smooth second-order analysis of the type we develop here, which is a powerful tool in variational analysis and nonsmooth optimization, is not well-known in the robust recovery literature, and are thus of particular importance in this context.

\item Finally, Section~\ref{sec:stronggroup} is devoted to a comprehensive treatment of unique/strong solutions for group-sparsity minimization problems \cite{RRN12,RF08,YL06}, an important class of \eqref{BP0} where $J$ is the $\ell_1/\ell_2$ norm. Sufficient conditions for solution uniqueness for group-sparsity minimization problems were achieved in \cite{CR13, G11, JKL15, RRN12, RF08, V14, VPF15}, but none of them are necessary. To find a characterization for solution uniqueness, we first obtain a closed form for the descent cone (Theorem~\ref{Uniq}). To the best of our knowledge, this was an open question in the literature, though the {\em statistical dimension} of this cone is more computable \cite{ALMT14}. By relying again on second-order analysis \cite{BS00, RW98}, we show that a unique solution to group-sparsity minimization problems \eqref{BP0} and \eqref{Las3} is indeed a strong solution (Theorem~\ref{BP12}). This result is interesting in two ways: first, the function in \eqref{BP0} is not piecewise linear-quadratic, which rules out the approach developed in \cite{BLN19} to unify solution uniqueness and strong minima; second, strong solutions are possibly more natural than sharp solutions for non-polyhedral minimization problems whenever solution uniqueness occurs. Moreover, we establish a quantitative condition for checking unique/strong solutions to group-sparsity problems that is equivalent to solving a smooth convex optimization  problem. This convex problem can be solved via available packages such as \texttt{cvxopt} as we will illustrate in the numerical results of Section~\ref{sec:numerics}. Consequently, we show that solution uniqueness to group-sparsity problem is equivalent to the robust recovery with rate $\mathcal{O}(\sqrt{\delta})$ for both problems \eqref{BP1} and \eqref{Las3}. 

\end{enumerate}

\begin{Remark}
Although we do not have a rigorous study for the case where $J$ is the nuclear norm in this paper, we provide several examples showing that a unique solution to \eqref{BP0} with the nuclear norm is neither sharp nor strong. This raises the challenge of understanding the interplay between unique, sharp, and strong solutions in this case. Some open questions will be discussed in Section~\ref{sec:conclusion} for further research in this direction.
\end{Remark}

%% file: tex/sec_prelim.tex
\section{Preliminaries}
\label{sec:prelim}
%%%%%%%%%%%%%%%%%%%%%%%%%%%%%%%%%%%%%%%%%%%%%%%%%%%%%%%%

Throughout the paper, $\R^n$ is the Euclidean space with dimension $n$. In $\R^n$ we denote the inner product by $\dotp{\cdot}{\cdot}$, the Euclidean norm by $\norm{\cdot}$, and the closed ball of radius $r > 0$ centered at $x \in \R^n$ by $\B_r(x)$, and $\Id$ the identity operator.

%%%%%%%%%%%%%%%%%%%%%%%%%%%%%%%%%%%%%%%%%%%%%%%%%%%%%%%%
\subsection{Some linear algebra} 
Given an $m\times n$ matrix $A$, $\Ker A$ (resp. $\Im A$) is the kernel/null (resp. image) space of $A$. The Moore-Penrose generalized inverse \cite[Page~423]{M00} of $A$ is denoted by $A^\dag$, which satisfies the property
\begin{equation}\label{MP1}
AA^\dag A=A.
\end{equation}
Moreover, we have 
\begin{equation}\label{cons}
\begin{gathered}
Ax=b \text{ is consistent if and only if } AA^\dag b=b, \qandq \\
\text{the set of solutions to $Ax=b$ is } A^\dag b+\Ker A.
\end{gathered}
\end{equation}
If $A$ is injective then $A^\dag=(A^\top A)^{-1}A^\top$, where $A^\top$ is the transpose of $A$. If $A$ is surjective then $A^\dag=A^\top (AA\top)^{-1}$.  $A^*$ will also denote the adjoint of $A$ as a linear operator.

Suppose that  $\ox$ is a solution to the consistent system $Ax=b$. The projection of $x\in \R^n$ onto the affine set $C \eqdef A^{-1}b$ is
\begin{equation}\label{Pro}
\proj_{C}(x)=(\Id-A^\dag A)(x-\ox)+\ox=x- A^\dag A(x-\ox);
\end{equation}
see, e.g., \cite[page~435--437]{M00}. Let $\|\cdot\|_\AA$ and $\|\cdot\|_\BB$ be two arbitrary norms in $\R^m$ and $\R^n$, respectively. The matrix operator norm $\|\cdot\|_{\AA,\BB}: \R^{m \times n} \to \R_+ $  is defined by
\begin{equation*}
\|A\|_{\AA,\BB} \eqdef  \max_{x\in \R^n\setminus\{0\} }\dfrac{\|Ax\|_\AA}{\|x\|_{\BB}}.
\end{equation*}
We write $\|A\|_{2,2}$ for $\|A\|$.  The Frobenius norm of $A$ is known as
\[
\|A\|_F=\sqrt{\Tr(A^*A)}.
\]
Recall that the nuclear norm of $A$ is 
\[
\|A\|_*=\sum_{k=1}^{\min(m,n)} \sigma_k,
\]
where $\sigma_k$, $1\le k\le \min(m,n)$ are all the singular values of $A$. 
%By the definition, for any $1\leq p, q\leq\infty$, $A\in\R^{m\times n}$ and $x\in \R^n$, one has
%\begin{equation}\label{norm_ineq}
%    \|Ax\|_p\leq \|A\|_{p,q}\|x\|_q.
%\end{equation}

%%%%%%%%%%%%%%%%%%%%%%%%%%%%%%%%%%%%%%%%%%%%%%%%%%%%%%%%
\subsection{Some convex analysis}
Let $\ph:\R^n \to \oR  \eqdef  \R\cup\{+\infty\}$ be a proper lower semi-continuous (lsc) convex function with domain $\dom \ph \eqdef \enscond{x\in \R^n}{\ph(x)<\infty}$. The subdifferential of $\ph$ at $x\in\dom \ph$ is defined by
\[
\partial \ph(x) \eqdef \enscond{v\in \R^n}{\ph(u) \geq \ph(x) + \dotp{v}{u-x}, \quad \forall u\in \R^n}. 
\]
When $\ph(\cdot)=\iota_C(\cdot)$, the indicator function of a nonempty closed convex set $C$, i.e., $\iota_C(x)=0$ if $x\in C$ and $+\infty$ otherwise, the subdifferential of $\iota_C$ at $x\in C$ is the {\em normal cone} to $C$ at $x$:
\begin{equation}\label{Nor}
N_C(x) \eqdef \enscond{v\in \R^n}{\dotp{v}{u-x} \le 0, \quad \forall u\in C}.
\end{equation}
For a set $\Omega \subset \R^n$, $\cone \Omega$ is its {\em conical hull}, $\Int \Omega$ is the interior of $\Omega$, and $\ri \Omega$ is the relative interior of the convex set $\Omega$. Given a convex cone $K \subset \R^n$, the {\em polar} of $K$ is defined to be the cone
\begin{equation}\label{Pola}
K^\circ \eqdef \enscond{v \in \R^n}{\dotp{v}{u} \le 0, \quad \forall u\in \Omega} .
\end{equation}
It is well known, see \cite[Corollary~6.21]{RW98} that $K^\circ$ is closed and convex, and that 
\begin{equation}\label{Dual}
K^\circ=\cl(K)^\circ \qandq K^{\circ\circ} = (K^\circ)^\circ= \cl (K),
\end{equation}
where $\cl$ denotes the topological {\em closure}. We will also use $\cl$ for the closure of a function, i.e., the topological closure of its epigraph.

The {\em gauge function} of $\Omega$ is  defined by 
\begin{equation}\label{ga}
\gg_\Omega(u) \eqdef \inf\enscond{r>0}{u\in rC}, \qforq u\in \R^n.
\end{equation}
The {\em support function} to $\Omega$ is denoted by
\begin{equation}\label{Sup}
\sigma_\Omega(v) \eqdef \sup_{x\in \Omega}\dotp{v}{x}, \qforq v\in \R^n.
\end{equation}
It is a proper lsc convex function as soon as $\Omega$ is non-empty.
The {\em polar set} of $\Omega$ is given by
\begin{equation}\label{Polar}
    \Omega^\circ \eqdef \enscond{v\in \R^n}{\sigma_ C(v) \le 1}.
\end{equation}
When $\Omega$ is a non-empty convex set, $\Omega^\circ$ is a non-empty closed convex set containing the origin. When $\Omega$ is also closed and contains the origin, it is well-known that
\begin{equation}\label{PoGa}
\gg_\Omega=\sigma_{\Omega^\circ},
\end{equation}
see, e.g., \cite[Corollary~15.1.2]{R70}. In this case, $\gamma_C$ is a non-negative lsc convex and positively homogenous function.

%%%%%%%%%%%%%%%%%%%%%%%%%%%%%%%%%%%%%%%%%%%%%%%%%%%%%%%%
\subsection{Sharp and strong minima} 
Let us first recall the definition of sharp and strong solutions in \cite{BS00,C78,P87,RW98} that are playing central roles throughout the paper.

\begin{Definition}[Sharp and strong minima] We say $\ox$ to be a sharp solution/minimizer to the (non necessarily convex) function $\ph$ with a constant $c>0$ if there exists $\ve>0$ such that 
\begin{equation}\label{Sha}
    \ph(x) \geq \ph(\ox) + c\norm{x-\ox},\quad \forall x \in  \B_\ve(\ox).
\end{equation}
Moreover, $\ox$ is said to be a strong solution/minimizer with a constant $\kk>0$ if  there exists $\delta>0$ such that 
\begin{equation}\label{Str}
    \ph(x) \geq \ph(\ox) + \frac{\kk}{2}\|x-\ox\|^2,\quad \forall x \in \B_\delta(\ox).
\end{equation}
\end{Definition}
Properties \eqref{Sha} and \eqref{Str} are also known as respectively, first order and second order (or quadratic) growth properties.

Strong minima is certainly weaker than sharp minima. Moreover, if $\ox$ is a sharp or strong solution to $\ph$, it is a unique local minimizer to $\ph$ (and the unique minimizer if $\ph$ is convex). When the function $\ph$ is a convex {\em piecewise linear convex}, i.e., its epigraph is a polyhedral,  $\ox$ is a unique solution to $\ph$ if and only if  it is a sharp solution. Furthermore, if $\ph$ is a convex {\em piecewise linear-quadratic} function in the sense that its domain is a union of finitely many polyhedral sets, relative to each of which $\ph$ is a quadratic function, $\ox$ is a strong solution if and only if  it is a unique solution, see, e.g., \cite{BLN19}.

To characterize sharp and strong solutions, it is typical to use directional derivative and second  subderivative; see, e.g., \cite{BS00,RW98}. %\JF{In the following, I've insisted on the fact that $\ph$ is convex as otherwise, there are many subtelties to be taken care of for the properties hereafter to hold true.}
%Actually, without convexity, this rather requires $\ph$ to be subdifferentially regular and strictly continuous/locally Lipschitz continuous (see \cite[Theorem~19.6]{RW98}). Think of $|\cdot|^p, p \in ]0,1[$, at the origin}
\begin{Definition}[directional derivative and second subderivative]\label{DD12}  
Let $\ph:\R^n \to \oR$ be a proper lsc convex function. The directional derivative of $\ph$ at $\ox \in \dom \ph$ is the function $d\ph(\ox):\R^n\to\oR$ defined by
\begin{equation}\label{sub1}
d\ph(\ox)(w) \eqdef \lim_{t\dn 0}\dfrac{\ph(\ox+tw)-\ph(\ox)}{t} \qforq w\in \R^n.
\end{equation}
The second subderivative of $\ph$ at $\ox \in \dom \ph$ for $\ov\in \partial \ph(\ox)$ is the function $d^2 \ph(\ox|\ov):\R^n\to \oR_+$  defined by   
\begin{equation}\label{Subd}
d^2\ph(\ox|\ov)(w) \eqdef \liminf_{t\dn 0, w^\prime \to w}\dfrac{\ph(\ox+tw^\prime)-\ph(\ox)-t\dotp{\ov}{ w^\prime}}{\frac{1}{2}t^2} \qforq w\in \R^n.
\end{equation}
\end{Definition}
It is well-known that $d\ph(\ox)(w)=\max\enscond{\dotp{v}{w}}{v\in \partial \ph(\ox)}$ if $\ph$ is continuous at $\ox$. The calculation of $d^2\ph(\ox|\ov)(w)$ is quite involved in general; see, e.g., \cite{BS00,RW98,MR11}. Note that from \cite[Proposition~13.5 and Proposition~13.20]{RW98}, the function $d^2\ph(\ox|\ov)$ is non-negative, lsc, positively homogeneous of degree 2, and 
\begin{equation}\label{dom2}
\dom d^2\ph(\ox|\ov) \subset \KK(\ox|\ov) \eqdef \enscond{w\in \R^n}{d\ph(\ox)(w)=\dotp{\ov}{w}}.
\end{equation}
If, moreover, $\ph$ is twice-differentiable at $\ox$, we have (see \cite[Example~13.8]{RW98})
\begin{equation}\label{smooth}
d^2\ph(\ox|\nabla \ph(\ox))(w)=\dotp{w}{\nabla^2\ph(\ox)w}.
\end{equation}
The following sum rule for second subderivatives will be useful. 
\begin{Lemma}
Let $\ph, \phi:\R^n\to\oR$ be proper lsc convex functions with $\ox\in \dom \ph\cap\dom \phi$. Suppose that $\ph$ is twice-differentiable at $\ox\in {\rm int}\,(\dom \ph)$ and $\ov\in \partial (\ph+\phi)(\ox)=\nabla \ph(\ox)+\partial \phi(\ox)$. Then we have
\begin{equation}\label{sum}
d^2(\ph+\phi)(\ox|\ov)(w)=\la w,\nabla^2\ph(\ox)w\ra+d^2\phi(\ox| \ov-\nabla \ph(\ox))(w) \qforallq w\in \R^n.
\end{equation}
\end{Lemma}

\begin{proof}
Use convexity of $\ph,\phi$ and twice-differentiability of $\ph$ at $\ox$ into \cite[Proposition~13.19]{RW98}.
\end{proof}

The next result taken from \cite[Lemma~3, Chapter~5]{P87} and \cite[Theorem~13.24]{RW98} characterize sharp and strong minima.

\begin{Lemma}[Characterization of sharp and strong minima]\label{Fa} 
Let $\ph:\R^n \to \oR$ be a proper lsc convex function with $\ox\in \dom \ph$. We have: 
\begin{enumerate}[label={\rm (\roman*)}]
\item \label{lemFa:item1}
$\ox$ is a sharp minimizer to $\ph$ if and only if there exists $c>0$ such that $d\ph(\ox)(w)\ge c\|w\|$ for all $w\in \R^n$, i.e., $d\ph(\ox)(w)>0$ for all $w\in \R^n\setminus\{0\}$. 
   % \item[{\rm (b)}] There exist $c,\delta>0$ such that ${\rm dist}\,(0,\ph(x))\ge c$ for all $x\in \B_\ve(\ox)\setminus\{\ox\}$ .

\item \label{lemFa:item2}
$\ox$ is a strong minimizer to $\ph$ if and only if $0\in \partial \ph(\ox)$ and there exists $\kk>0$ such that $d^2\ph(\ox|0)(w)\ge \kk\|w\|^2$ for all $w\in \R^n$, which means
\[
\Ker d^2\ph(\ox|0) \eqdef \enscond{w\in\R^n}{d^2\ph(\ox|0)(w)=0}=\{0\}.
\]
\end{enumerate}
\end{Lemma}

It is important to observe that although the sharpness of a minimizer as defined in \eqref{Sha} is a local property, it is actually a global one for convex problems. 
\begin{Lemma}[Global sharp minima]\label{Glo} 
Suppose that $\ox$ is a sharp solution to the proper lsc convex function $\ph$. Then there exists $c>0$ such that 
\[
\ph(x)-\ph(\ox)\ge c\|x-\ox\|\quad \qforallq x\in \R^n.
\]
\end{Lemma}

\begin{proof} 
Since $\ox$ is a sharp solution to $\ph$, we have from Lemma~\ref{Fa}\ref{lemFa:item1} that $\exists c > 0$ such that $d\ph(\ox)(w)\ge c\|w\|$, $\forall  w\in \R^n$. Thus, by the characterization of the directional derivative for convex functions, we have for any $x\in \R^n$
\[
\ph(x)-\ph(\ox)=\ph(\ox+(x-\ox))-\ph(\ox)\ge \inf_{t>0}\dfrac{\ph(\ox+t(x-\ox))-\ph(\ox)}{t}=d\ph(\ox)(x-\ox)\ge c\|x-\ox\|. 
\]
The proof is complete.
\end{proof}

%% file: tex/sec_unique.tex
%%%%%%%%%%%%%%%%%%%%%%%%%%%%%%%%%%%%%%%%%%%%%%%%%%%%%%%%
\section{Unique, sharp, and strong solutions for robust recovery}
\label{sec:unique}
%%%%%%%%%%%%%%%%%%%%%%%%%%%%%%%%%%%%%%%%%%%%%%%%%%%%%%%%

%%%%%%%%%%%%%%%%%%%%%%%%%%%%%%%%%%%%%%%%%%%%%%%%%%%%%%%%
\subsection{Uniqueness and robust recovery}
Let us consider the optimization problem in \eqref{BP0}
% \begin{equation*}\label{ABP}
%     \min\quad J(x)\quad \mbox{subject to}\quad \Phi x=\Phi x_0,
% \end{equation*}
where we suppose that $J:\R^n \to \oR$ is a non-negative lsc function but not necessarily convex. We denote by $J_\infty: \R^n \to \oR$ the {\em asymptotic (or horizon) function} \cite{AT03} associated with $J$, which is defined by
\begin{equation}\label{AF}
J_\infty(w) \eqdef \liminf_{w^\prime \to w, t\to \infty}\dfrac{J(tw^\prime)}{t} .
\end{equation}
Throughout this section, we assume that $J$ satisfies 
\begin{equation}\label{AC}
\Ker J_\infty\cap \Ker \Phi=\{0\},
\end{equation}
where $\Ker J_\infty \eqdef \enscond{w\in \R^n}{J_\infty(w)=0}$, and the range of $J_\infty$ is on $\R_+$ since $J$ is non-negative. It then follows from \cite[Corollary~3.1.2]{AT03} that \eqref{AC} ensures that problem \eqref{BP0} has a non-empty compact set of minimizers. 

Define $\Psi:\R^n\to \oR$ as
\begin{equation}\label{Psi}
\Psi(x) \eqdef  J(x)+\iota_{\Phi^{-1}(\Phi x_0)}(x)\quad \mbox{for all}\quad x\in \R^n.
\end{equation}
We say that $x_0$ is a unique, sharp, or strong solution to \eqref{BP0} if it is a unique, sharp, or strong solution to the function $\Psi$, respectively.

Our first result shows that solution uniqueness is sufficient for robust recovery. If $J$ is also convex, then uniqueness is also necessary. The sufficiency part of our result generalizes the corresponding part of \cite[Theorem~3.5]{HKPS07} to nonconvex problems. 

\begin{Proposition}[Solution uniqueness for robust recovery]\label{NRC} 
Suppose that $J$ is a non-negative lsc function which satisfies \eqref{AC}. 
\begin{enumerate}[label={\rm (\roman*)}]
\item If $x_0$ is the unique solution to problem \eqref{BP0} then: \label{NRC:claim1}
\begin{enumerate}[label={\rm (\alph*)}]
\item any solution $x_\delta$ to problem \eqref{BP1} with $\|y-y_0\| \le \delta$ converges to $x_0$ as $\delta \to 0$, \label{NRC:claim1a}
    
\item for any constant $c_1>0$, any solution $x_{\mu}$ to problem \eqref{Las3} with $\mu=c_1\delta$ and $\|y-y_0\|\le \delta$ converges to $x_0$ as $\delta \to 0$. \label{NRC:claim1b}
\end{enumerate}
\item Conversely, if $J$ is convex and $\dom J = \R^n$, then: \label{NRC:claim2}
\begin{enumerate}[label={\rm (\alph*)}] 
\item \ref{NRC:claim1}-\ref{NRC:claim1a} implies $x_0$ is the unique solution to \eqref{BP0}. 
\item \ref{NRC:claim1}-\ref{NRC:claim1b} implies $x_0$ is the unique solution to \eqref{BP0}. 
\end{enumerate}
\end{enumerate}
\end{Proposition}
\begin{proof} 
\begin{enumerate}[label={\rm (\roman*)}]
\item Suppose that $x_0$ is the unique solution to problem \eqref{BP0}. To justify \ref{NRC:claim1}-\ref{NRC:claim1a}, define the function $g \eqdef \iota_{\B_\delta(y)} \circ \Phi$. We have, using \cite[Propositions~2.6.1, 2.6.3 and 2.1.2]{AT03}, that
\begin{eqnarray*}
(J+g)_\infty(w) \geq J_\infty(w)+g_\infty(w) = J_\infty(w)+(\iota_{\B_\delta(y)})_\infty(\Phi w) 
&=& J_\infty(w)+\iota_{\{0\}}(\Phi w) \\
&=& J_\infty(w)+\iota_{\Ker \Phi}(w) \geq 0 .
\end{eqnarray*}
Thus, under \eqref{AC}, the set of solutions to \eqref{BP1} is nonempty and compact thanks to \cite[Corollary~3.1.2]{AT03}. If $\{x_\delta\}_{\delta>0}$ is unbounded, without passing to subsequences suppose that $\|x_\delta\|\to \infty$ and  $\dfrac{x_\delta}{\|x_\delta\|} \to w$ as $\delta \dn 0$ with $\|w\|=1$. We have
\[ 
0 \le \anorm{\Phi \pa{\frac{x_\delta}{\|x_\delta\|}}-\frac{y}{\|x_\delta\|}} \le \frac{\delta}{\|x_\delta\|},
\] 
and after passing to the limit, we get that $w\in \Ker \Phi$. Moreover, since $x_0$ is a feasible point of \eqref{BP1}, we have $J(x_\delta) \le J(x_0)$, and thus
\[
0=\liminf_{\delta\dn 0}\dfrac{J(x_0)}{\|x_\delta\|}\ge \liminf_{\delta\dn 0}\dfrac{J(x_\delta)}{\|x_\delta\|}\ge J_\infty(w)\ge 0, 
\]
which means that $w \in \Ker J_\infty$. This entails that  $0 \neq w \in \Ker J_\infty\cap \Ker \Phi$, which  contradicts \eqref{AC}. Thus $\{x_\delta\}_{\delta>0}$ is bounded. Pick any subsequence $\{x_{\delta_k}\}_{k \in \N}$ of $\{x_\delta\}_{\delta>0}$ converging to some $\bar{x}$ as $\delta_k \dn 0$. From lower semicontinuity of $J$ and of the norm, we get
\[
J(\bar{x}) \leq \liminf_{k} J(x_{\delta_k}) \leq J(x_0) \qandq 0 \leq \norm{y_0 - \Phi\bar{x}} \leq \liminf_{k} \norm{y_0 - \Phi x_{\delta_k}} \leq 2 \lim_{k} \delta_k = 0 .
\]
This means that $\bar{x}$ is a solution of \eqref{BP0}, and by uniqueness of the minimizer, $\bar{x}=x_0$. Since the subsequence $\{x_{\delta_k}\}_{k \in \N}$ was arbitrary, we have $x_\delta \to x_0$ as $\delta\dn 0$, whence claim \ref{NRC:claim1}-\ref{NRC:claim1a} follows.\\

The proof of \ref{NRC:claim1}-\ref{NRC:claim1b} uses a similar reasoning. Let $g = \frac{1}{2}\norm{y-\cdot}^2 \circ \Phi$. One has, using \cite[Propositions~2.6.1, 2.6.3 and 2.6.7]{AT03}, that
\begin{eqnarray*}
(\mu J+g)_\infty(w) \geq \mu J_\infty(w)+g_\infty(w) = \mu J_\infty(w)+\frac{1}{2}\pa{\norm{\cdot}^2}_\infty(\Phi w) 
&=& \mu J_\infty(w)+\iota_{\{0\}}(\Phi w) \\
&=& \mu J_\infty(w)+\iota_{\Ker \Phi}(w) \geq 0 .
\end{eqnarray*}
Hence the set of optimal solutions to problem \eqref{Las3} is nonempty and compact according to \cite[Corollary~3.1.2]{AT03}. Similarly to part \ref{NRC:claim1}-\ref{NRC:claim1a}, we claim that  $\{x_\mu\}_{\mu > 0}$ is bounded. By contradiction, suppose again that $\|x_\mu\|\to \infty$  and $ \dfrac{x_\mu}{\|x_\mu\|} \to z$ with $\|z\|=1$. Note that by optimality of $x_\mu$
\begin{equation}\label{Sim}
\mu J(x_\mu) \leq \frac{1}{2}\|\Phi x_\mu-y\|^2+\mu J(x_\mu) \leq \frac{1}{2}\|\Phi x_0-y\|^2+\mu J(x_0) \leq \frac{\delta^2}{2}+\mu J(x_0) .
\end{equation}
As $\mu=c_1\delta$, we obtain from \eqref{Sim} that 
\[
0=\liminf_{\delta\dn 0}\dfrac{J(x_0)+\frac{\delta}{2c_1}}{\|x_\mu\|}\ge \liminf_{\delta\dn 0}\dfrac{J(x_\mu)}{\|x_\mu\|}\ge J_\infty(z)\ge 0, 
\]
which means $z\in \Ker J_\infty$. Furthermore, \eqref{Sim} also entails 
\[
0\le \frac{1}{2}\anorm{\Phi\pa{\frac{x_\mu}{\|x_\mu\|}}-\frac{y}{\|x_\mu\|}}^2 \le \frac{1}{\|x_\mu\|^2 } \pa{\frac{\delta^2}{2}+\mu J(x_0)},
\]
which tells us, after passing to the limit, that $0 \neq z \in \Ker J_\infty \cap \Ker \Phi $, hence contradicting \eqref{AC}. Thus $\{x_\mu\}_{\mu > 0}$ is bounded. Arguing as in \ref{NRC:claim1}-\ref{NRC:claim1a}, we use lower semicontinuity of $J$ and of the norm to show that any cluster point of $\{x_\mu\}_{\mu > 0}$ is a solution of \eqref{BP0}, and deduce claim \ref{NRC:claim1}-\ref{NRC:claim1b} thanks to uniqueness of the minimizer $x_0$.
%Let $\{x_{\mu_q}\}$ be an arbitrary subsequence of $\{x_\mu\}$ converging to some $v_0$. Taking $\delta\to 0$, we derive from \eqref{Sim} that $J(x_0)\ge J(v_0)$ and $\Phi v_0=y_0$. Since $x_0$ is the unique solution to problem \eqref{BP0}, $x_0=v_0$. As the  subsequence $\{x_{\mu_q}\}$ is arbitrary,  $\{x_\mu\}$ converges to $x_0$ as $\delta\dn 0$. This verifies \ref{NRC:claim2}. 

\item Suppose now that $J$ is also convex and $\dom J = \R^n$ (hence continuous). If \ref{NRC:claim1}-\ref{NRC:claim1a} holds, we have 
\begin{equation}\label{BP4}
J(x_\delta) \leq J(x) \qforanyq x \qstq \|y - \Phi x\| \le \delta \qandq \|y-y_0\|\le \delta.
\end{equation}
It follows that 
\[
J(x_\delta) \leq J(x) \qforanyq x \qstq \Phi x=y_0.
\]
By letting $\delta\dn 0$, and using continuity of $J$, we obtain that $x_0$ is a minimizer to \eqref{BP0}. Let $\ox$ be an arbitrary minimizer to problem \eqref{BP0}. As $\dom J=\R^n$, Fermat's rule for \eqref{BP0} and subdifferential calculus gives
\begin{equation}\label{eq:optcondBP0}
0\in \partial \Psi(\ox)=\partial J(\ox)+N_{\Phi^{-1}y_0}(\ox)=\partial J(\ox)+\Im \Phi^*.
\end{equation}
Or, equivalently, there exists a dual multiplier $\eta\in \R^m$ such that $\Phi^*\eta\in \partial J(\ox)$\footnote{This is known as the Source Condition \cite{scherzer2009variational}.}. If $\eta=0$, we have $0\in  \partial J(\ox)$, which means that $\ox$ is a minimizer to $J$. Thus, $\ox$ is also a minimizer to problem \eqref{BP1} for any $\delta>0$. By \ref{NRC:claim1}-\ref{NRC:claim1a}, $\ox-x_0\to 0$ as $\delta\dn 0$, i.e., $\ox=x_0$. 

If $\eta\neq 0$, take any $\delta>0$ and define $\mu \eqdef \dfrac{\delta}{\|\eta\|}$ and $y \eqdef y_0+\mu \eta=\Phi \ox+\mu \eta$. It follows that $\|y - \Phi\ox\| = \|y - y_0\| = \delta$ and that
\begin{equation}\label{eq:fermatox}
\frac{1}{\mu}\Phi^*(y - \Phi\ox)= \Phi^*\eta\in \partial J(\ox).
\end{equation}
On the other hand, under our assumption, we have that $x^\star$ is a solution to \eqref{BP1} if and only if
\[
0 \in \partial J(x^\star)  + \Phi^* N_{\B_\delta(y)}(\Phi x^\star) ,
\]
where $N_{\B_\delta(y)}(\Phi x^\star) = \R_+(\Phi x^\star-y)$ whenever $\norm{y-\Phi x^\star}=\delta$. This is precisely the optimality condition verified by $\ox$ in \eqref{eq:fermatox}, which shows that $\ox$ is also a solution to \eqref{BP1}. By \ref{NRC:claim1}-\ref{NRC:claim1a} again, we have $\ox=x_0$. Since the choice of  $\ox$ is arbitrary, $x_0$ is the unique solution to \eqref{BP0}.

Finally, suppose that \ref{NRC:claim1}-\ref{NRC:claim1b} is satisfied. We have for any $x \in \Phi^{-1}y_0$ that
\[
\mu J(x_\mu) \leq \frac{1}{2}\|\Phi x-y\|^2+\mu J(x) \leq \frac{\delta^2}{2}+\mu J(x) .
\]
Dividing both sides by $\mu=c_1\delta$ and letting $\delta\to 0$ entails $J(x_0) \leq J(x)$ for any  $x\in \Phi^{-1}y_0$, which is equivalent to saying that $x_0$ is a solution to problem \eqref{BP0}. The proof of uniqueness follows similar lines as for problem \eqref{BP1} with the choice $y = \Phi \ox + \mu \eta$ and $\mu=c_1\delta$. Note that under our assumption, Fermat's rule for problem \eqref{Las3} is \eqref{eq:fermatox} without further qualification condition unlike what is required for \eqref{BP1}.
\end{enumerate}
\end{proof} 

\begin{Remark}\label{rem:coeruniqueness}
In the convex case, the condition \eqref{AC} is superfluous in Proposition~\ref{NRC}\ref{NRC:claim1}. Indeed, \eqref{AC} is equivalent in this case to the fact the set of minimizers of \eqref{BP0} is a non-empty compact set, which is obviously the case when $x_0$ is the unique minimizer.
\end{Remark}

%%%%%%%%%%%%%%%%%%%%%%%%%%%%%%%%%%%%%%%%%%%%%%%%%%%%%%%%
\subsection{Sharp minima and robust recovery}

A sufficient and necessary condition for uniqueness to problem \eqref{BP0}, together with its implication for robust recovery, is studied in \cite{CRPW12} via the key geometric notion of {\em descent cone} that we recall now.
\begin{Definition}[Descent cone] The descent cone of the function $J$ at $x_0$ is defined by 
\begin{equation}\label{Des}
\DD_J(x_0) \eqdef \cone\enscond{x-x_0}{J(x) \le J(x_0)} .
\end{equation}
\end{Definition}
 
The following proposition provides a characterization for solution uniqueness whose proof can be found in \cite[Proposition~2.1]{CRPW12}.

\begin{Proposition}[Descent cone for solution uniqueness \cite{ALMT14,CRPW12}]\label{DeUn}  Then
$x_0$ is the unique solution to \eqref{BP0} if and only if $\Ker \Phi \cap \DD_J(x_0)=\{0\}$. 
\end{Proposition}

In \cite[Proposition~2.2]{CRPW12} (see also \cite[Proposition~2.6]{Tropp15}), the authors also prove a robust recovery result for solutions of \eqref{BP1} via the descent cone. %Since our problem \eqref{BP0} is slightly different  from the model in \cite{CRPW12}, we provide a short proof.

\begin{Proposition}[Robust recovery via descent cone \cite{CRPW12}] \label{Pop} Let  
Suppose that there exists some $\al>0$ such that $\|\Phi w\|\ge \al\|w\|$ for all $w\in \DD_J(x_0)$. Then any solution $x_\delta$ to problem \eqref{BP1} satisfies
\begin{equation}\label{Est7}
\|x_\delta-x_0\|\le  \frac{2\delta}{\al}.
\end{equation}
\end{Proposition}

The parameter $\alpha$ is known as the minimum conic singular value of $\Phi$ with respect to the cone $\DD_J(x_0)$. One may observe that the condition in Proposition~\ref{Pop} is not the same with the one for solution uniqueness in Proposition~\ref{DeUn}. The reason is that the descent cone $\DD_J(x_0)$ is not necessarily closed in general. To see this, we first unveil the relation between the descent and {\em critical} cones. Denote 
\begin{equation*}
L_{x_0}(J) \eqdef \enscond{x\in \R^n}{J(x)\le J(x_0)}
\end{equation*}
the {\em sublevel set} of $J$ at $x_0$. The {\em critical cone} of a convex function $J$ at $x_0$ is 
\begin{equation}\label{D0}
\CC_J(x_0) \eqdef \enscond{w \in \R^n}{dJ(x_0)(w)\le 0} .
\end{equation} 

\begin{Lemma}[Relationship between descent cone and critical cone]\label{lem:tancritcones}
Suppose that $J$ is a continuous convex function. Then,
\begin{equation}\label{CDinc}
\DD_J(x_0) \subset \cl (\DD_J(x_0)) \subset \CC_J(x_0).
\end{equation}
If moreover $0 \notin \partial J(x_0)$, then
\begin{equation}\label{CD4}
\cl (\DD_J(x_0)) = \CC_J(x_0).
\end{equation}
\end{Lemma}
\begin{proof}
We have from \eqref{Dual} that
\begin{equation}\label{ClDinc}
\cl (\DD_J(x_0))=(\DD_J(x_0))^{\circ\circ}=\pa{(\cl \DD_J(x_0))^{\circ}}^\circ=(N_{L_{x_0}(J)}(x_0))^\circ \subset \CC_J(x_0) ,
\end{equation}
where we used \eqref{Dual} in the first two equalities, polarity in the third, and convexity and the first claim of \cite[Proposition~10.3]{RW98} in the last equality. 
If in addition $0 \notin \partial J(x_0)$, then using \cite[Corollary~23.7.1]{R70} and the last claim of \cite[Proposition~10.3]{RW98}, the inclusion in \eqref{ClDinc} becomes
\begin{equation}\label{ClD}
\cl (\DD_J(x_0))=(N_{L_{x_0}(J)}(x_0))^\circ=(\cone \partial J(x_0))^\circ=\CC_J(x_0) .
\end{equation}
\end{proof}

% \begin{proof} 
% Suppose that $x_0$ is the unique solution to \eqref{BP0} but $0 \neq w \in \Ker \Phi \cap \DD_\AA(x_0)$. Hence there exists some $t>0$ such that $J(x_0+tw) \le J(x_0)$ and $\Phi(x_0+tw)=\Phi x_0$, or equivalently that $x_0+tw$ is a solution \eqref{BP0}, contradicting uniquness of $x_0$. Conversely, suppose that $x_0$ is not the unique solution but $\Ker \Phi \cap \DD_\AA(x_0) = \{0\}$. We know then that there exists $\ox \neq x_0$ such that $J(\ox) \leq J(x_0)$ and $\ox - x_0 \in \Ker \Phi$. This implies that $\ox-x_0\in \Ker \Phi \cap \DD_\AA(x_0) \neq \{0\}$, leading again to a contradiction. The proof is complete. 
% \end{proof} 

When $J$ is also positively homogenous, the requirement $0 \notin \partial J(x_0)$ in Lemma~\ref{lem:tancritcones} can be dropped.

\begin{Proposition}\label{prop:tancritcones}
Suppose that $J=\sigma_{C}$, where $C \subset \R^n$ is a non-empty compact convex set with $0 \in \ri C$. Then, \eqref{CD4} holds at any $x_0 \in \R^n$.
\end{Proposition}
\begin{proof} 
Thanks to the assumptions on $C$, we have from \cite[Theorem~13.2 and Corollary~13.3.1]{R70} that $J$ convex, positively and homogenous and finite everywhere, hence continuous. We now only consider the case where $0 \in \partial J(x_0)$ as, otherwise, the result claim follows from Lemma~\ref{lem:tancritcones}. By \cite[Lemma~1]{VGFP15}, $J$ is also non-negative and $\Ker J$ is a linear subspace. This entails that $0 \in \partial J(x_0)$ is equivalent to $x_0 \in \Ker J$, and thus from \eqref{Des}, that 
\[
\DD_J(x_0) = \cone (\Ker J - x_0) = \Ker J .
\]
On the other hand, we have from \eqref{D0} that
\[
\CC_J(x_0) = (\partial J(x_0))^\circ = (\cone \partial J(x_0))^\circ=N_{\partial J(x_0)}(0) .
\]
Moreover, by \cite[Corollary~8.25]{RW98} 
\[
\partial J(x_0) = C \cap \{x_0\}^\perp .
\]
Combining this with \cite[Theorem~6.42]{RW98} (recall that $C$ is convex and $0 \in \ri C \cap \{x_0\}^\perp$), we get
\[
\CC_J(x_0) = N_{C \cap \{x_0\}^\perp}(0) = N_C(0) + \R x_0 = \Ker J + \R x_0 = \Ker J ,
\]
where the last equality comes from the fact that $\Ker J$ is a linear subspace containing $x_0$. This completes the proof. 
\end{proof}

We now turn to showing that the condition in Proposition~\ref{Pop} actually means that $x_0$ is a sharp solution to \eqref{BP0}. 

\begin{Proposition}[Characterization of sharp minima via the descent and critical cones]\label{Shar} 
Let $J$ be a continuous convex function.
\begin{enumerate}[label={\rm (\roman*)}]
\item $x_0$ is a sharp solution to \eqref{BP0} if and only if $\Ker \Phi\cap \CC_J(x_0)=\{0\}$. \label{Shar:claim1}
\item If $J=\sigma_{C}$, where $C \subset \R^n$ is a non-empty compact convex set with $0 \in \ri C$, then $x_0$ is a sharp solution to \eqref{BP0} if and only if there exists $\al>0$ such that $\|\Phi w\|\ge \al\|w\|$ for all $w\in \DD_J(x_0)$. \label{Shar:claim2}
\end{enumerate}

\end{Proposition}
\begin{proof} 
\begin{enumerate}[label={\rm (\roman*)}]
\item Recall $\Psi$ from \eqref{Psi}. We have 
\begin{equation}\label{dP}
d\Psi (x_0)(w)=dJ(x_0)(w)+\iota_{\Ker \Phi}(w).
\end{equation}
By Lemma~\ref{Fa}, $x_0$ is a sharp solution to \eqref{BP0} if and only if $d\Psi (x_0)(w)>0$ for all $w\in \R^n\setminus\{0\}$. Combining this with \eqref{D0} and \eqref{dP} leads the claimed equivalence. 

\item Assume $x_0$ is a sharp solution to \eqref{BP0}. Since $\CC_J(x_0)$ is closed, we have from claim~\ref{Shar:claim1} that there exists $\al>0$ such that $\|\Phi w\|\ge \al\|w\|$ for all $w\in \CC_J(x_0)$. We then get the first implication from the inclusion \eqref{CDinc}. Conversely, if there is $\al>0$ such that $\|\Phi w\|\geq \al\|w\|$ for all $w\in \DD_J(x_0)$, we have $\|\Phi w\|\ge \al\|w\|$ for all $w \in \CC_J(x_0)$ thanks to Proposition~\ref{prop:tancritcones}. It follows that $\Ker \Phi \cap \CC_J(x_0)=\{0\}$. This implies, in view of claim~\ref{Shar:claim1}, that $x_0$ is a sharp solution to \eqref{BP0}. 
\end{enumerate}
\end{proof}

The main lesson we learn from the above, in particular Proposition~\ref{DeUn} and Proposition~\ref{Shar}, is that there is a gap between solution uniqueness and solution sharpness (hence robust recovery via Proposition~\ref{Pop}) for problem \eqref{BP0}. This difference lies in the closedness of the descent cone. Indeed, when the descent cone $\DD_J(x_0)$ is closed, it coincides with the critical cone $\CC_J(x_0)$. Unfortunately, in general, $\DD_J(x_0)$ is not closed. A prominent example where $\DD_J(x_0)$ may fail to be closed is that of the $\ell_1/\ell_2$ norm very popular to promote group sparsity. The following example is a prelude of our precise formula for the descent cone of the $\ell_1/\ell_2$ norm in Theorem~\ref{Uniq}, in which the interior of the critical cone $\CC_J(x_0)$ plays a significant role.

\begin{Example}[Gap between solution uniqueness and sharpness for group-sparsity]\label{Ex12}{\rm  Consider the following $\ell_1/\ell_2$ norm minimization problem:
\begin{equation}\label{Ex123}
\min_{x\in \R^3} J(x)=\sqrt{x_1^2+x_2^2}+|x_3| \qsubjq \Phi x=y_0 ,
\end{equation}
with $\Phi=\begin{pmatrix}1& 1 &0\\ 1& 0 &-1\end{pmatrix}$, $x_0=(0,1,0)^\top$, and $y_0=\Phi x_0=(1,0)^\top$. We have $\Ker \Phi = \R(1,-1,1)^\top$, and feasible points of \eqref{Ex123} then take the form $x=x_0+t(1,-1,1)^\top$ for any $t \in \R$. For all such points, we have 
\begin{eqnarray*}
J(x)-J(x_0)=\sqrt{t^2+(1-t)^2}+|t|-1=\dfrac{t^2+2(|t|-t)}{\sqrt{t^2+(1-t)^2}+1-|t|},
\end{eqnarray*}
which tells us that $x_0$ is the unique solution to \eqref{Ex12}. In fact, simple calculation shows that $x_0$ is a strong but not a sharp solution to \eqref{Ex12}. For $J$ in \eqref{Ex123}, Figure~\ref{fig:conesl1l2} illustrates the difference between the descent cone computed by Theorem~\ref{Uniq} and the critical cone \eqref{D0} whose computation is straightforward in this case. It is obvious that $\Ker \Phi\cap \CC_J(x_0)\neq \{0\}$, though $\Ker \Phi\cap \DD_J(x_0)= \{0\}$ (which is indeed equivalent to solution uniqueness). 
\begin{figure}[!ht]
\centering
\includegraphics[width=0.5\linewidth]{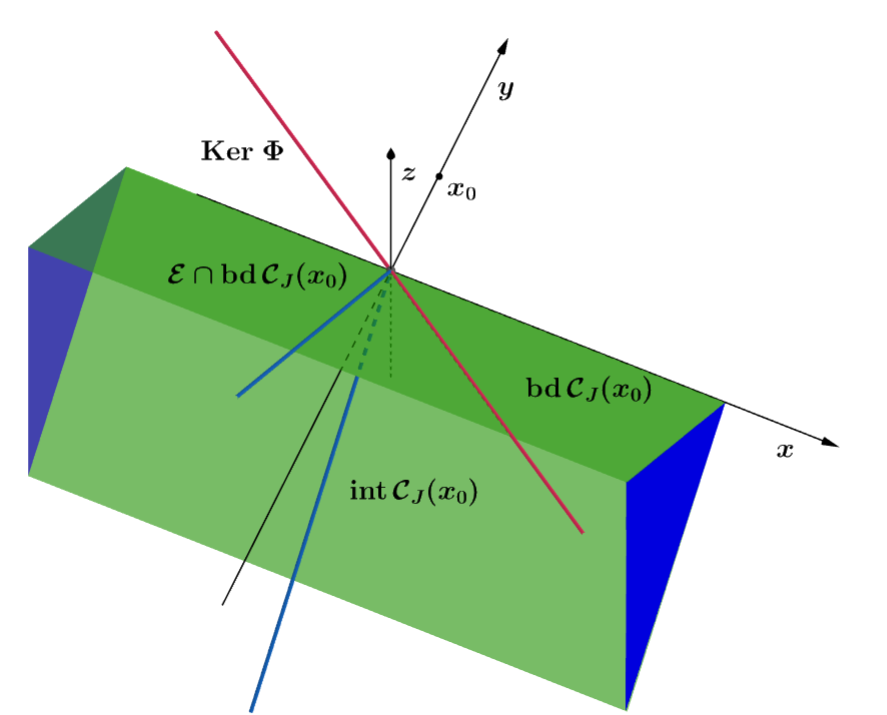}
\caption{The descent cone $\DD_J (x_0)$ and critical cone $\CC_J (x_0)$ for $J$ in \eqref{Ex123}. Here, $\DD_J (x_0)=(\EE\cap\bd \CC_J(x_0))\cup\Int \CC_J(x_0)$, where $\EE=\enscond{(0,x_2,x_3)}{x_2, x_3 \in\R}$ (see \eqref{C}) and $\CC_J(x_0)=\enscond{(x_1,x_2,x_3)}{x_2+|x_3|\leq 0}$. The green shaded surface is $\bd \CC_J(x_0)$ and the blue half-lines correspond to $\EE \cap \bd \CC_J(x_0)$.}
\label{fig:conesl1l2}
 \end{figure}
} \qed
\end{Example}

%\[
%\lim_{x_1\downarrow 0}\dfrac{J(x)-J(x_0)}{\|x-x_0\|}=\lim_{x_1\downarrow 0}\dfrac{x_1}{\sqrt{3}(\sqrt{x_1^2+(1-x_1)^2}+1-x_1)}=0.
%\]

Inspired by Proposition~\ref{Pop}, we show next that sharpness of the minimizer guarantees robust recovery for both problems \eqref{BP1} and \eqref{Las3} with linear rate. Unlike \cite[Proposition~2.2]{CRPW12} and many other results \cite{FPVDS13, GHS11, V14, VPF15, ZYY16} in this direction, we do not need convexity of $J$. In fact, the key is to assume that $x_0$ is a unique and sharp solution to problem \eqref{BP0}.

\begin{Theorem}[Sharp minima for robust recovery]\label{Rob} 
Suppose that $J$ is a non-negative lsc function which satisfies \eqref{AC}, and is Lipschitz continuous around $x_0$ with modulus $L$. If  $x_0$ is the unique and sharp solution to \eqref{BP0} with constant $c>0$, then for all $\delta > 0$ sufficiently small, the following statements hold:
\begin{enumerate}[label={\rm (\roman*)}]
\item \label{Rob:claim1}
Any solution $x_\delta$ to problem \eqref{BP1} satisfies
\begin{equation}\label{Est2}
\|x_\delta-x_0\|\le \frac{2(L+c)\|\Phi^\dag\|}{c}\delta.
\end{equation}   
\item \label{Rob:claim2}
For any $c_1>0$ and $\mu=c_1\delta$, any minimizer $x_{\mu}$ to \eqref{Las3} satisfies
\begin{equation}\label{Est}
\|x_\mu-x_0\|\le \frac{c_1}{2c} \left(\frac{1}{c_1}+(c+L)\|\Phi^\dag\|\right)^2\delta.
\end{equation}
\end{enumerate}
\end{Theorem}

\begin{proof}
\begin{enumerate}[label={\rm (\roman*)}]
\item $L$-Lipschitz continuity of $J$ around $x_0$ entails that $\exists \ve > 0$ such that 
\begin{equation}\label{Lip}
|J(x)-J(z)|\le L\|x-z\| \qforallq x,z\in \B_\ve(x_0). 
\end{equation}
Since $x_0$ is a sharp solution to \eqref{BP0} with constant $c$, there exists $\eta>0$ such that 
\begin{equation}\label{InSharp}
J(x)-J(x_0)\ge c\|x-x_0\| \qforallq x\in \Phi^{-1}y_0\cap\B_\eta(x_0). 
\end{equation} 
Let $\ox_\delta$ be the projection of $x_\delta$ onto $\Phi^{-1}y_0$. Thus
\begin{equation}\label{eq:normprojxdel}
\norm{\ox_\delta-x_0} \leq \norm{x_\delta-\ox_\delta} + \norm{x_\delta-x_0} \leq 2 \norm{x_\delta-x_0} .
\end{equation} 
Since $x_0$ is the unique solution to problem \eqref{BP0}, Proposition~\ref{NRC}\ref{NRC:claim1} tells us that for all $\delta$ small enough, $x_\delta \in  \B_{\min\{\ve,\eta\}/2}(x_0)$ and, by \eqref{eq:normprojxdel}, $\ox_\delta \in  \B_{\min\{\ve,\eta\}}(x_0)$, i.e., $\ox_\delta, x_\delta\in \B_{\min\{\ve,\eta\}}(x_0)$ with $\ve, \eta$ in \eqref{Lip} and \eqref{InSharp}. We then have
\begin{align*}
L\|\ox_\delta-x_\delta\|
&\ge J(\ox_\delta)-J(x_\delta) 	& \text{by \eqref{Lip}}, \\
&\ge J(\ox_\delta)-J(x_0)	& \text{by optimality of $x_0$}, \\
&\ge c\|\ox_\delta-x_0\| 	& \text{by \eqref{InSharp}}, \\
&\ge c\|x_\delta-x_0\|-c\|x_\delta-\ox_\delta\|.
\end{align*}
Combining this with the projection formula \eqref{Pro} tells us that 
\begin{align*}
c\|x_\delta-x_0\|
&\le (L+c)\|x_\delta-\ox_\delta\|=(L+c)\|\Phi^\dag \Phi (x_\delta -x_0)\|\\
&\le (L+c)\|\Phi^\dag\|\|\Phi (x_\delta -x_0)\|\\
&\le (L+c)\|\Phi^\dag\|(\|y-\Phi x_\delta\|+\|y-\Phi x_0\|)\\
&\le (L+c)\|\Phi^\dag\|2\delta,
\end{align*}
which is \eqref{Est2}.

\item Optimality of $x_\mu$ gives 
\begin{align}%\label{Sim1}
J(x_\mu)-J(x_0) 
&\leq \frac{1}{2\mu}\pa{\|y-\Phi x_0\|^2 - \|y - \Phi x_\mu\|^2} \nonumber \\
&\leq \frac{1}{2\mu}\pa{\|y-\Phi x_0\|^2 - \pa{\|y-\Phi x_0\| - \|\Phi (x_\mu-x_0)\|}^2} \nonumber \\
&= \frac{1}{2\mu} \pa{2\delta\|\Phi (x_\mu-x_0)\| - \|\Phi (x_\mu-x_0)\|^2} \label{Low} .
\end{align}
%which tells us that $\|D^*x_\mu\|\le  J(x)+\frac{\|y-y_0\|^2}{\mu}\le  J(x)+\frac{\delta}{c_1}$ for any $x\in \Phi^{-1}\oy$. Take any subsequence  $\{x_{\mu_n}\}$ such that  $D^*x_{\mu_q}$ is converging to some $D^*u_0$. The latter inequality also implies that $\|D^*u_0\|\le  J(x)$ for any $x\in \Phi^{-1}y_0$. Hence $u_0$ is a solution to problem \eqref{BP0}. Since $x_0$ is a sharp solution to \eqref{BP0}, we have $u_0=x_0$. Note that the choice of converging subsequence $\{x_{\mu_q}\}$ is arbitrary, $D^* x_\mu \to D^*x_0$ for sufficiently small $\delta$. Since $D$ is surjective, $x_\mu \to x_0$ as $\delta\dn 0$. Without loss of generality, we suppose $\|x_\mu-x_0\|\le \ve$.

Let $\ox_\mu$ be the projection of $x_\mu$ to $\Phi^{-1}y_0$. Inequality \eqref{eq:normprojxdel} remains valid for $x_\mu$ and $\ox_\mu$ replacing $x_\delta$ and $\ox_\delta$ respectively. Moreover, arguing as in the proof of claim~\ref{Rob:claim1}, but now invoking Proposition~\ref{NRC}\ref{NRC:claim2}, we infer that for $\delta$ sufficiently small, $x_\mu,\ox_\mu\in \B_{\min\{\ve,\eta\}}(\ox)$, where $\ve,\eta$ are those in \eqref{Lip} and \eqref{InSharp}. Denote for short $\al=\|\Phi(x_\mu-x_0)\|$. We then get
\begin{align*}
J(x_\mu)-J(x_0)
&=J(\ox_\mu)-J(x_0)+J(x_\mu)-J(\ox_\mu)\\
&\ge c\|\ox_\mu-x_0\|-L\|x_\mu-\ox_\mu\| 		& \text{by \eqref{InSharp} and \eqref{Lip}}, \\
&\ge c\|x_\mu-x_0\|-(c+L)\|x_\mu-\ox_\mu\|		& \\
&=   c\|x_\mu-x_0\|-(c+L)\|\Phi^\dag \Phi(x_\mu-x_0)\| 	& \text{by \eqref{Pro}}, \\
&\ge c\|x_\mu-x_0\|-(c+L)\|\Phi^\dag\|\al .
\end{align*}
This together with \eqref{Low} gives us that
\begin{equation}\label{beta}
\begin{aligned}
c\|x_\mu-x_0\|
&\leq \frac{2\delta\al-\al^2}{2\mu}+(c+L)\|\Phi^\dag\| \al\\
&=\pa{\frac{\delta}{\mu}+(c+L)\|\Phi^\dag\|}\al-\frac{\al^2}{2\mu}\\
&\disp\le\frac{\mu}{2} \pa{\frac{\delta}{\mu}+(c+L)\|\Phi^\dag\|}^2,
\end{aligned}
\end{equation}
which is clearly \eqref{Est} for the choice $\mu=c_1\delta$.
\end{enumerate}
\end{proof} 

The bounds \eqref{Est2} and \eqref{Est} in Theorem~\ref{Rob} hold without convexity with the proviso that $\delta$ is sufficiently small. However, when the function $J$ is convex and continuous, the bounds are satisfied for any $\delta>0$.

\begin{Corollary}[Sharp minima for robust recovery under convexity]\label{RobCon} 
Suppose that $J$ is a convex continuous function. If $x_0$ is a sharp solution to \eqref{BP0}, then statements \ref{Rob:claim1} and \ref{Rob:claim2} of Theorem~\ref{Rob} hold for any $\delta>0$ with some constant $L>0$.
\end{Corollary}
\begin{proof} 
When $J$ is a continuous convex function and $x_0$ is a sharp solution to \eqref{BP0}, it follows from Lemma~\ref{Glo} that the sharpness property \eqref{InSharp} is global and $x_0$ is also the unique solution to \eqref{BP0}. In turn, $J$ satisfies \eqref{AC} in view of Remark~\ref{rem:coeruniqueness}. Therefore, the nets $\{x_\delta\}_{\delta > 0},\{x_\mu\}_{\mu > 0}$ are bounded as proved in Proposition~\ref{NRC}. Thus, there exists some $R>0$ such that $x_0,x_\delta,\ox_\delta, x_\mu,\ox_\mu\in \B_R(0)$. Convexity and continuity of $J$ also imply that it is Lipschitz continuous on $\B_{2R}(x_0)$ with some Lipschitz constant $L>0$. Overall, this tells us that we can take $\ve, \eta$ in \eqref{Lip} and \eqref{InSharp} in the proof of Theorem~\ref{Rob} by $2R$. The rest of the proof remains valid, whence our claim follows. 
\end{proof}

%\begin{Remark}\label{Relax}
%%%%%%%%%%%%%%%%%%%%%%%%%%%%%%%%%%%%%%
\paragraph{Discussion of related work}
It is worth noting that Corollary~\ref{RobCon} covers many results in \cite{FPVDS13, GHS11, V14, DT14, VPF15, ZYY16}. When $J=\|\cdot\|_\AA$ is a norm in $\R^n$, claim \ref{Rob:claim1} of Corollary~\ref{RobCon} is exactly \cite[Proposition~2.2]{CRPW12} (see Proposition~\ref{Pop}) thanks to Proposition~\ref{Shar}\ref{Shar:claim2}. In this case the Lipschitz constant $L$ of $J$ is $\|\Id\|_{\AA,2}$.  
 
For the case  $J=\|\cdot\|_{\ell_1}$,  Corollary~\ref{RobCon} returns \cite[Theorem~4.7]{GHS11} (see also \cite[Proposition~1]{DT14}) whose proof is less transparent and involves deriving linear convergence rate of the Bregman divergence of $J$, together with the characterizations of solution uniqueness to $\ell_1$ optimization problems via the so-called restricted injectivity and non-degenerate source condition.

When $J(x)=\|D^*x\|_\AA$ with $D$ being an $n\times p$ matrix and $\|\cdot\|_\AA$ being a norm in $\R^p$, the Lipschitz constant $L$ of $J$ is $\|\Id\|_{\AA,2}\|D^*\|$. Corollary~\ref{RobCon} covers \cite[Theorem~2]{FPVDS13}. Another result in this direction is \cite[Theorem~2]{ZYY16}, which only obtains linear rate for $\|D^*(x_\mu-x_0)\|$ with extra nontrivial assumptions on $D$. 

When $J$ is a general convex continuous regularizer as in Corollary~\ref{RobCon}, \cite{V14,VPF15} use the so-called  restricted injectivity and non-degenerate source condition at $x_0$ to obtain robust recovery. In our forthcoming Theorem~\ref{mtheo}, we will show that these two conditions are equivalent to $x_0$ being a sharp solution\footnote{The authors in \cite{V14,VPF15} have already proved that these two conditions are sufficient for $x_0$ to be a sharp solution.}. It means that Corollary~\ref{RobCon} is equivalent to \cite[Theorem~6.1]{V14} or \cite[Theorem~2]{VPF15}. However, our path to proving robust recovery is radically different. On the one hand, \cite[Theorem~6.1]{V14} generalizes the proof strategy initiated in \cite{GHS11} and use decomposability of $J$ and other arguments that heavily rely on convexity of $J$. On the other hand, Corollary~\ref{RobCon} provides a direct proof that involve natural geometrical properties of $J$ around $x_0$. Most importantly, Theorem~\ref{Rob} suggests that robust recovery with linear rate occurs  without convexity and reveals the crucial role played by sharpness of the minimizer to achieve robust recovery with linear rate. Comparing the constants in our bounds \eqref{Est}-\eqref{Est2} and those in \cite{V14,VPF15}, those in \cite{V14,VPF15} depend for instance on a dual certificate and its "distance" to degeneracy, while ours depend on the sharpness parameter $c$ which is not trivial to characterize at first glance. In  Theorem~\ref{mtheo}, we will provide an estimation for $c$ via the so-called {\em Source Identity} that can be computed numerically.

%The approach is employed later in Theorem~\ref{SRob} when sharp minima is replaced by strong minima.   

%Although the Lipschitz modulus $L$ could be estimated for any continuous convex functions \cite[Proposition~8.28]{BC11}, it is rather complicated. However, the similar bound to \eqref{Est} in \cite[Theorem~6.1]{V14}

%Although descent cones can be used to prove (i) as in Proposition~\ref{Shar}, it is not clear how to use it to prove (ii). Our technique in Theorem~\ref{Rob} is uniform for both (i) and (ii) by using the projection to the affine constraint. The same idea will be  employed in Theorem~\ref{SRob}.

%Another important thing about Theorem~\ref{Rob} is that both bounds obtained in \eqref{Est2} and \eqref{Est} are explicit. The only constant that is not trivial is $c$. In  Theorem~\ref{mtheo} we provide an estimation for this number via the so-called {\em Source Identity} that can be computed numerically. The bounds in the aforementioned results  \cite{FPVDS13, GHS11, V14, VPF15, ZYY16} have  hidden constants  that are hard to estimate. \end{proof} 

%\end{Remark}

%%%%%%%%%%%%%%%%%%%%%%%%%%%%%%%%%%%%%%%%%%%%%%%%%%%%%%%%
\subsection{Strong minima and robust recovery}
A natural question is to whether robust recovery is still possible if the sharp minima property is replaced by the weaker strong minima property. We here show that the answer is affirmative, but at the price of a slower rate of convergence. 

\begin{Theorem}[Strong minima for robust recovery]\label{SRob} 
Suppose that $J$ is a non-negative lsc function which satisfies \eqref{AC}, and is Lipschitz continuous around $x_0$ with modulus $L$. If  $x_0$ is the unique and strong solution to \eqref{BP0} with constant $\kk>0$, then for all $\delta > 0$ sufficiently small, the following statements hold:
\begin{enumerate}[label={\rm (\roman*)}]
\item \label{SRob:claim1}
Any solution $x_\delta$ to problem \eqref{BP1} satisfies
\begin{eqnarray}\label{Est4}
\|x_\delta-x_0\| \le 2\pa{\frac{1}{\kk}L\|\Phi^\dag\| \delta+\|\Phi^\dag\|^2\delta^2}^{\frac{1}{2}}.
\end{eqnarray}
    
\item \label{SRob:claim2}
For any $c_1>0$ and $\mu=c_1\delta$, any minimizer $x_{\mu}$ to \eqref{Las3} obeys
\begin{equation}\label{Est3}    \|x_\mu-x_0\|\le\sqrt{\dfrac{c_1}{(1-c_1\kk\|\Phi^\dag\|^2\delta)\kk}}\left(\dfrac{1}{c_1}+L\|\Phi^\dag\|\right)\delta^\frac{1}{2}.
\end{equation}
\end{enumerate}
\end{Theorem}
\begin{proof} 
\begin{enumerate}[label={\rm (\roman*)}]
\item Since $x_0$ is a strong solution to \eqref{BP0} with constant $\kk>0$, there exists $\nu>0$ such that 
\begin{equation}\label{Stron}
J(x)-J(x_0)\ge \frac{\kk}{2}\|x-x_0\|^2 \qforallq x\in \Phi^{-1}y_0\cap\B_\nu(x_0). 
\end{equation} 
Let $\ox_\delta$ be the projection of $x_\delta$ onto $\Phi^{-1}y_0$. Since $x_0$ is the unique solution to problem \eqref{BP0}, we argue as in the proof of Theorem~\ref{Rob}\ref{Rob:claim1} to show that for all $\delta$ small enough, one has $\ox_\delta, x_\delta \in \B_{\min\{\ve,\nu\}}(x_0)$ with $\ve, \nu$ in \eqref{Lip} and \eqref{Stron}.
%Since $x_0$ is the unique solution to problem \eqref{BP0}, we have $\Ker D^*\cap \Ker \Phi=\{0\}$. Hence there exists $\ell>0$ such that 
%\begin{equation}\label{DPhi}
  %   J(x)\ge \ell\|x\|\quad \mbox{for all}\quad x\in \Ker \Phi.
%\end{equation}
%To justify (i), let $\ox_\delta$ be the projection of $x_\delta$ to $\Phi^{-1}y_0$. 
%Note that $J(x_0)\ge J(x_\delta)$. Since $J$ is coercive, the sequence $\{x_\delta\}$ is bounded. Pick any subsequence $\{x_{\delta_q}\}$ of $\{x_\delta\}$ converging to some $u_0$ as $q\to \infty$. We get that $J(x_0)\ge J(u_0)$ and $\Phi u_0=y_0$. Since $x_0$ is the unique solution to \eqref{BP0}, we have $u_0=x_0$. 
%\begin{eqnarray}\label{Phid2}\begin{array}{ll}
%J(x_0)&\disp\ge J(x_\delta)\ge J(\ox_\delta)-L\|x_\delta-\ox_\delta\|\\
%&\ge \disp J(\ox_\delta)- L\|\Phi^\dag\Phi(x_\delta-x_0)\|\quad \mbox{(by \eqref{Prox})}\\
%&\ge \disp J(\ox_\delta)- L\|\Phi^\dag\|\cdot\|\Phi(x_\delta-x_0)\|\\
%&\ge\disp J(\ox_\delta)- L\|\Phi^\dag\|2\delta.
%\end{array}
%\end{eqnarray}
%Moreover, we obtain from \eqref{DPhi} that
%\begin{equation}\label{Bo}
%J(\ox_\delta)\ge \|D^*(\ox_\delta-x_0)\|_\AA-J(x_0)\ge \ell\|\ox_\delta-x_0\|-J(x_0).
%\end{equation}
We then have 
\begin{align*}
L\|\ox_\delta-x_\delta\|
&\ge J(\ox_\delta)-J(x_\delta) 	& \text{by \eqref{Lip}}, \\
&\ge J(\ox_\delta)-J(x_0)	& \text{by optimality of $x_0$}, \\
&\ge \frac{\kk}{2}\|\ox_\delta-x_0\|^2 	& \text{by \eqref{Stron}}, \\
&= \frac{\kk}{2}\pa{\|x_\delta-x_0\|^2-\|\ox_\delta-x_\delta\|^2}. & \text{by \eqref{Pro}}.
\end{align*}
This together with \eqref{Pro} again tells us that 
\begin{align*}
\frac{\kk}{2}\|x_\delta-x_0\|^2
&\le  L \|\ox_\delta-x_\delta\|+\frac{\kk}{2}\|\ox_\delta-x_\delta\|^2\\
&= L\|\Phi^\dag \Phi (x_\delta -x_0)\|+\frac{\kk}{2}\|\Phi^\dag \Phi (x_\delta -x_0)\|^2\\
&\le L\|\Phi^\dag\|\|\Phi x_\delta-y_0\|+\frac{\kk}{2}\|\Phi^\dag\|^2\|\Phi x_\delta-y_0\|^2\\
&\le L\|\Phi^\dag\| 2\delta+\frac{\kk}{2}\|\Phi^\dag\|^24\delta^2,
\end{align*} 
which is \eqref{Est4}.

\item Let $\ox_\mu$ be the projection of $x_\mu$ to $\Phi^{-1}y_0$. Arguing as above, we infer that for $\delta$ sufficiently small, $x_\mu,\ox_\mu\in \B_{\min\{\ve,\nu\}}(\ox)$, where $\ve,\eta$ are those in \eqref{Lip} and \eqref{Stron}. Set for short $\al=\|\Phi(x_\mu-x_0)\|$. We then have
\begin{align*}
J(x_\mu)-J(x_0)
&=J(\ox_\mu)-J(x_0)+J(x_\mu)-J(\ox_\mu)\\
&\ge \frac{\kk}{2}\|\ox_\mu-x_0\|^2-L\|\ox_\mu-x_\mu\|			& \text{by \eqref{Stron} and \eqref{Lip}}, \\
&=\frac{\kk}{2}\pa{\|x_\mu-x_0\|^2-\|\ox_\mu-x_\mu\|^2} - L\|\ox_\mu-x_\mu\|	& \text{by \eqref{Pro}}, \\
&=\frac{\kk}{2}\|x_\mu-x_0\|^2 - \|\Phi^\dag\Phi(x_\mu-x_0)\|^2 - L\|\Phi^\dag\Phi(x_\mu-x_0)\| & \text{by \eqref{Pro}}, \\
&\geq\frac{\kk}{2}\pa{\|x_\mu-x_0\|^2 - \|\Phi^\dag\|^2\alpha^2} - L\|\Phi^\dag\|\alpha . &
\end{align*}
Combining this with \eqref{Low}, we arrive at
\[
\frac{\kk}{2}\pa{\|x_\mu-x_0\|^2-\|\Phi^\dag\|^2\al^2}-L\|\Phi^\dag\|\al \leq \frac{2\al\delta-\al^2}{2\mu} .
\]
Making $\delta$ smaller if necessary so that $\mu = c_1\delta < \pa{\kk\norm{\Phi^\dag}^2}^{-1}$, we get
\begin{align*}
\frac{\kk}{2}\|x_\mu-x_0\|^2
&\le \pa{\frac{\delta}{\mu}+L\|\Phi^\dag\|}\al-\pa{\frac{1}{2\mu}-\frac{\kk}{2}\|\Phi^\dag\|^2}\al^2\\
&\le \frac{1}{2}\pa{\frac{1}{\mu}-\kk\|\Phi^\dag\|^2}^{-1}\pa{\frac{\delta}{\mu}+L\|\Phi^\dag\|}^2,
\end{align*}
which is \eqref{Est4} after some simple algebra.
\end{enumerate}
\end{proof}

When the regularizer $J$ is convex, the Lipschitz property of J around $x_0$ is equivalent to the continuity of $J$ at $x_0$; see, e.g., \cite[Theorem~8.29]{BC11}. Moreover, the assumption that $x_0$ is the unique minimizer of \eqref{BP0} holds trivially when $x_0$ is a strong minimizer. We then obtain the following corollary of Theorem~\ref{SRob}.
 
\begin{Corollary}\label{SRobCon}
Suppose that $J$ is a convex function that is continuous at $x_0$. If $x_0$ is a strong solution to problem \eqref{BP0} then, for $\delta$ sufficiently small, the convergence rate $\mathcal{O}(\delta^\frac{1}{2})$ of Theorem~\ref{SRob} holds. 
\end{Corollary}

\begin{Remark}
Observe that unlike sharp minima, even for the convex case, the property of strong minima is only local. This is the reason we still require $\delta$ to be small enough in Corollary~\ref{SRobCon}. This is in stark contrast to Corollary~\ref{RobCon} where sharpness is globalized under convexity.
\end{Remark}

\begin{Remark}
A natural and open question is whether the rate $\mathcal{O}(\delta^\frac{1}{2})$ in Theorem~\ref{SRob} is optimal under the strong solution property. Though we do not have a clear answer yet, we believe that his may depend on the type of regularizer. For the case of Euclidean norm $J(x)=\|x\|$, we could prove that the rate could be improved to $\mathcal{O}(\delta)$. We omit the details here for the sake of brevity.
\end{Remark}

%% file: tex/sec_sharp.tex
\section{Nondegeneracy, restricted injectivity and sharp minima}
\label{sec:sharp}
%%%%%%%%%%%%%%%%%%%%%%%%%%%%%%%%%%%%%%%%%%%%%%%%%%%%%%%%

As shown in Proposition~\ref{Shar}, sharp minima can be characterized via the descent and critical cones. However, it is not easy to verify them numerically, especially when the dimension is large\footnote{For some operators $\Phi$ drawn from some appropriate random ensembles, one can show, using , that sample complexity bounds in \cite{CRPW12,Tropp15} are sufficient for sharpness to hold with high probability.}. In this section, we mainly derive quantitative characterizations for  sharp minima to problem \eqref{BP0}. 

Throughout this section, we suppose that $J$ takes the analysis-type form
\begin{equation}\label{eq:Janalysis}
J = J_0 \circ D^* ,
\end{equation}
where $D^*: \R^n \to \R^p$ is linear, and $J_0$ is a non-negative, convex and continuous function. The reason we take this form is twofold. First, this is in preparation for Section~\ref{sec:stronggroup} to make the presentation easier there. Second, though we could derive the decomposability properties (see shortly) of $J$ from those of $J_0$ using the framework in \cite{VGFP15}, our analysis in this section will involve a dual multiplier which is not the same as the one in that paper.

For some linear subspace $V$ of $\R^n$, we will use the shorhand notation $w_V = \proj_V w$ for $w \in \R^n$ and $D_V = D \proj_V$ for a linear operator $D$.

%%%%%%%%%%%%%%%%%%%%%%%%%%%
\subsection{Subdifferential decomposability}
The following definition taken from \cite[Definition~3]{VGFP15} is instrumental in our study. 
\begin{Definition}[Model Tangent Subspace]\label{MTan} 
Denote the {\em model vector} $e$ of $J_0$ at $u_0 \eqdef D^* x_0$ is the projection of $0$ onto the affine hull of $\partial J_0(u_0)$
\begin{equation}\label{e}
e \eqdef \proj_{\aff (\partial J_0(u_0))}(0).
\end{equation}
The {\em model tangent subspace} $T$ at $u_0$ associated to $J_0$ is defined by
\begin{equation}\label{TS}
T \eqdef S^\perp \qwhereq S \eqdef  \aff (\partial J_0(u_0))-e.
\end{equation}
\end{Definition}
Obviously, $e \in T$.

\begin{Lemma}[Decomposability, {\cite[Theorem~1]{VGFP15}}]\label{Decom} 
Let $u_0\in \R^n\setminus\{0\}$ and $v_0$ be a vector in $\ri \partial J_0({u_0})$. Then 
\begin{equation}\label{Decop}
\partial J_0(u_0)=\enscond{v\in \R^n}{v_T=e, \gg_{C}(\proj_S(v-v_0))\le 1},
\end{equation}
where $C \eqdef  \partial J_0(u_0)-v_0$ and $\gg_C$ is the gauge function of $C$ defined in \eqref{ga}. Moreover $v\in \ri\partial J_0(u_0)$ if and only if $v_T=e$ and $\gg_C(\proj_S(v-v_0))<1$. 
\end{Lemma}

%For $x_0\in \R^n$, we have $\ri (\partial J(x_0))\neq \emptyset$. Pick a vector $v_0\in \ri (\partial J(x_0))$

%As discussed in Introduction part, solution uniqueness to \eqref{BP0} plays an essential role in robust and exact recovery. It could be characterized by using the so-called {\em descent tangent cone} \cite{CRPW12}. Although descent tangent cone is an important notion that is used widely in many papers, it does not have a closed form. More explicit characterizations are desired. In the case of $\ell_1$ problems, solution uniqueness to \eqref{BP0} is characterized via the Restricted Injectivity and Source Condition in \cite{GHS11}. But in $\ell_1$ problem a solution is unique if and only if it is a sharp solution  imply sharp minima. This is the lying reason that many results on robust and exact recovery only need solution uniqueness. For different frameworks, unique solution is distinct from sharp solution. Let us consider the following simple example, which tells that a unique solution to a problem is not a sharp solution, but a strong solution. Indeed, for the case of group sparsity minimization problem, solution uniqueness and strong minima are equivalent;  see further study in our Section~5.

We then have the following necessary and sufficient optimality condition for problem \eqref{BP0}.
\begin{Lemma}\label{lem:optcondBP0analysis}
$x_0$ is an optimal solution to problem \eqref{BP0} with $J$ in \eqref{eq:Janalysis} if and only if there exists $v\in \R^p$ such that 
\begin{equation}\label{eq:optcondBP0analysis}
v_T=e, \gg_C(v_S-\proj_Sv_0)\le 1, \qandq Dv\in \Im\Phi^* .
\end{equation}
\end{Lemma}
\begin{proof}
Combine \eqref{eq:optcondBP0} and Lemma~\ref{Decom}.
\end{proof}

%%%%%%%%%%%%%%%%%%%%%%%%%%%
\subsection{Quantitative characterization of sharp minima}
We start by providing a quantitative condition for checking optimality. 

\begin{Proposition}[Quantitative characterization of optimality]\label{Sol}  
$x_0$ is an optimal solution to problem \eqref{BP0} with $J$ in \eqref{eq:Janalysis} if and only if the following two conditions are satisfied
\begin{eqnarray}
&&ND_{S}(ND_{S})^\dag ND_Te=ND_Te \qandq  \label{yz}\\
&&\rho(e) \eqdef \min_{u\in \Ker ND_{S}}\gg_C(-(ND_{S})^\dag ND_Te-\proj_Sv_0+u)\le 1,  \label{xyz}
\end{eqnarray}
where $N$ is a matrix satisfying $\Ker N=\Im \Phi^*$ and $\Ker N^*=\{0\}$. The minimum in \eqref{xyz} is well-defined and achieved at some $u \in \Ker ND_{S} \cap S$.
\end{Proposition}
\begin{proof} 
Suppose that $x_0$ is a solution to \eqref{BP0}. It follows from Lemma~\ref{lem:optcondBP0analysis} that there exists $v\in \R^p$ such that 
\[
0=NDv=N(D_T e + D_{S}v_S) \Longrightarrow ND_{S}v_S=-ND_Te .
\]
This is equivalent, via \eqref{cons}, to 
\[
ND_{S}(ND_{S})^\dag ND_Te=ND_Te \qandq v_S \in -(ND_{S})^\dag ND_Te + \Ker ND_{S}.
\]
This tells us, using again Lemma~\ref{lem:optcondBP0analysis}, that
\[
\rho(e)\le \gg_{C}(v_S-\proj_Sv_0)\le 1,
\]
which proves the necessary part.

Conversely, suppose that both \eqref{yz} and \eqref{xyz} hold. Let $u^\star \in \Ker ND_{S}$ be a minimizer to \eqref{xyz}. Since $\dom \gg_C = S$ and $\gg_C$ is coercive on $S$ by \cite[Propopsition~2]{VGFP15}, $u^\star$ exists and belongs to $S$. Define $\ou \eqdef  -(ND_{S})^\dag NDe + u^\star$. This vector verifies
\begin{equation}\label{eq:ubargauge}
\gg_C(\ou-\proj_Sv_0)=\rho(e)\le1.
\end{equation}
Since $0 \in \ri C \subset C$, and in view of \eqref{eq:ubargauge}, we infer that $\ou-\proj_Sv_0\in C\subset S$. This implies $\ou\in S$. Note further from \eqref{yz} that
\[
ND(e+\ou)=ND(e+\proj_S\ou)=ND_Te-ND_{S}(ND_{S})^\dag ND_Te=0,
\]
which means $D(e+\ou)\in \Im \Phi^*$. Altogether, the vector $e+\ou$ verifies the properties \eqref{eq:optcondBP0analysis} in Lemma~\ref{lem:optcondBP0analysis} whence we deduce that $x_0$ is an optimal solution to problem \eqref{BP0}. The proof is complete. 
\end{proof}

Equation \eqref{yz} means that the linear system
\[
ND_{S}v=-ND_Te
\]
is consistent. When $D=\Id$, Lemma~\ref{lem:optcondBP0analysis} and Proposition~\ref{Sol} tell us that conditions \eqref{yz}-\eqref{xyz} are equivalent to the so-called Source Condition well-known in inverse problems (see \cite{scherzer2009variational,GHS11,V14,VPF15} and referencs therein). We call $\rho(e)$ the {\em Source Coefficient}. Under condition \eqref{yz}, $\rho(e)$ is indeed the optimal value to the following problem
\begin{equation*}\label{ro1}
\min_{v \in \R^p} \gg_C(v-\proj_Sv_0) \qsubjq ND_Sv=-ND_Te,
\end{equation*}
which is equivalent to the following convex optimization problem
\begin{equation}\label{ro2}
\min_{v \in \R^p}  \gg_C(v-\proj_Sv_0) \qsubjq  NDv=-ND_Te \qandq v\in S. 
\end{equation}

\begin{Remark}[Computing the Source Coefficient $\rho(e)$ for decomposable norms]\label{Deco}
One important class of regularizers is $J(x)=\|D^*x\|_\AA$, where $J_0(\cdot)=\|\cdot\|_\AA$ is a {\em decomposable} norm in $\R^p$. The class of decomposable norms is introduced in \cite{NRWY09,CR13} and was generalized in \cite{VGFP15} (coined strong gauges there). This class includes the $\ell_1$ norm, the $\ell_1/\ell_2$ norm, the nuclear norm but not the $\ell_\infty$ norm. According to \cite{CR13}, a norm $\|\cdot\|_\AA$ is called decomposable at $u_0 \eqdef D^*x_0$ if there is a subspace $V \subset \R^p$ and a vector $e_0 \in V$ such that
\begin{equation}\label{subA}
\partial \|u_0\|_{\AA}=\enscond{v\in \R^p}{v_V=e_0 \qandq \|v_{V^\perp}\|^*_{\AA}\le 1},
\end{equation}
where $\|\cdot\|^*_\AA$ is the dual norm to $\|\cdot\|_\AA$. From \cite[Proposition~7]{VGFP15}, it follows that \eqref{subA} complies with Lemma~\ref{Decom} by taking $V=T=S^\perp=\pa{\aff (\partial J_0(u_0))-e}^\perp$, $v_0 = e \in \ri \partial\|u_0\|_\AA$, in which case $C=\partial\|u_0\|_\AA-v_0=\enscond{v \in S}{\|v\|^*_\AA \le 1}$ and thus for any $v\in S$ 
\[
\gg_C(v)= \|v\|^*_\AA.
\]
Hence, problem \eqref{ro2} simplifies to
\begin{equation}\label{rhoz}
\min_{v \in \R^p} \|v\|^*_\AA \qsubjq NDv=-ND_Te \qandq v \in S.
\end{equation}
The matrix $N$ can be chosen from the singular value decomposition of $\Phi=U\Sigma V^*$ as $N=V_G^*$, where $V_G$ is the submatrix of $V$ whose columns are indexed by $G=\{r+1, r+2, \ldots, n\}$ with $r={\rm rank}\, \Phi$. This idea is slightly inspired from \cite{NDEG13,ZYY16} for the case where $\|\cdot\|_{\AA}$ is the $\ell_1$ norm. We will return in  Section~\ref{sec:stronggroup} to discussing further the use of the Source Coefficient to classify sharp and strong/unique solutions for the group-sparsity regularization. 
\end{Remark}

We are now in position to state the main result of this section, which provides equivalent characterizations of sharp minima. 

\begin{Theorem}[Characterizations of sharp solutions]\label{mtheo} 
Consider $J$ in \eqref{eq:Janalysis}. The following statements  are equivalent:
\begin{enumerate}[label={\rm (\roman*)}]
\item \label{mtheo:claim1}
$x_0$ is a sharp  solution to problem \eqref{BP0}. 

\item \label{mtheo:claim2}
The following {\em Restricted Injectivity} holds at $x_0$ 
\begin{equation}\label{GIR}
\Ker \Phi\cap \Ker D^*_S=\{0\}
\end{equation}  
and the {\em Source Coefficient} $\rho(e) < 1$.  

\item \label{mtheo:claim3}
The Restricted Injectivity \eqref{GIR} holds and the {\em Nondegenerate Source Condition} is satisfied at $x_0$ in the sense that 
\begin{equation}\label{NSC}
\exists v\in \R^p \qstq Dv\in \Im \Phi^* \qandq v\in \ri \partial J_0(x_0).
\end{equation}  
\end{enumerate}

Moreover, the sharpness constant $c$ at $x_0$ can be as large as
\begin{equation}\label{ShaCo}
c=(1-\rho(e))sc_1 > 0 \qwithq  c_1 \eqdef \min_{w\in \Ker \Phi\cap \mathbb{S}^{n-1}}\|D^*_{S}w\|>0
\end{equation}
and $s>0$ satisfying $\B_s(0)\cap S\subset C$. 
\end{Theorem}
\begin{proof} 
~[\ref{mtheo:claim1}$\Rightarrow$\ref{mtheo:claim2}]. Suppose that $x_0$ is a sharp solution to \eqref{BP0}. From Lemma~\ref{Fa}\ref{lemFa:item1}, this is equivalent to the existence of some $c>0$ such that $d\Psi(x_0)(w)\ge c\|w\|$ for all $w\in \R^n$. From \eqref{eq:optcondBP0} and Lemma~\ref{Decom}, we have
\begin{equation}\label{dPhi}
\begin{aligned}
    d\Psi(x_0)(w)
    &=dJ(x_0)(w)+\iota_{\Ker \Phi}(w)\\
    &=dJ_0(D^*x_0)(D^*w)+\iota_{\Ker \Phi}(w)\\
    &=\sigma_{\partial J_0(D^*x_0)}(D^*w)+\iota_{\Ker \Phi}(w)\\
    &=\sup\enscond{\dotp{v}{D^*w}}{v_T=e, v_S \in \proj_Sv_0+ C}+\iota_{\Ker \Phi}(w)\\
    &=\dotp{e}{D_T^*w} + \dotp{\proj_Sv_0}{D_S^*w} + \sup_{z \in C}\dotp{z}{D_S^*w}+\iota_{\Ker \Phi}(w) \\
    &=\dotp{e}{D_T^*w} + \dotp{\proj_Sv_0}{D_S^*w} + \sigma_C(D_S^*w)+\iota_{\Ker \Phi}(w) .
\end{aligned}
\end{equation}
For any $w=N^*u\in \Ker\Phi$, we obtain from \eqref{dPhi}, the inequality $d\Psi(x_0)(w)\ge c\|w\|$, and \eqref{yz} that  
\begin{equation}\label{Tes}
\begin{aligned}
    \sigma_C(D^*_Sw)-c\|w\|
    &\ge -\dotp{e}{D^*_TN^*u}-\dotp{\proj_Sv_0}{D^*_SN^*u}\\
    &=-\dotp{ND_Te}{u}-\dotp{\proj_Sv_0}{D^*_SN^*u}\\
    &=-\dotp{ND_S(ND_S)^\dag ND_Te}{u}-\dotp{ \proj_Sv_0}{D^*_SN^*u} \\
    &=-\dotp{(ND_S)^\dag ND_Te + \proj_Sv_0}{D^*_Sw} .
\end{aligned}
\end{equation}
If $w \in \Ker\Phi \cap \Ker D^*_S$, we deduce from \eqref{Tes} that $-c\|w\| \ge 0$, which means $w=0$. Thus the Restricted Injectivity \eqref{GIR} is satisfied. 

Moreover, as $J_0$ is a continuous convex function, $\partial J_0(D^* x_0)$ is compact, and so is $C$, whence it follows that 
\[
\sigma_C(D_S^*w) \le r\|D^*_Sw\|\le r\|D^*_S\|\|w\| \qwithq r \eqdef \max\enscond{\|z\|}{z \in C}.
\]
Combining the latter with \eqref{Tes} tells us that for any $v \in \Im D^*_SN^*$,
\begin{equation}\label{gau}
-\dotp{(ND_S)^\dag ND_Te+ \proj_Sv_0}{v} \leq \pa{1-\frac{c}{r\|D^*_S\|}}\sigma_C(v) .
\end{equation}

Let $C^\circ$ be the polar set of $C$ as defined in \eqref{Polar}, and set $K \eqdef C^\circ\cap S$. $K$ is a non-empty closed convex set. Since $0\in \ri C$, there exists $s>0$ such that and $\B_s(0)\cap S\subset C$. For any $v\in K\setminus\{0\}$, we have $s\frac{v}{\|v\|}\in \B_s(0)\cap S$ and thus, using \eqref{Polar}, $\dotp{v}{s\frac{v}{\|v\|}} \leq \sigma_{\B_s(0)\cap S}(v) \leq \sigma_C(v)\le 1$ since $v \in C^\circ$. It follows that $\|v\|\le \frac{1}{s}$. Hence, $K$ is a compact set. 

Let us bound from below the right hand side of \eqref{gau}. First, we have by Fenchel-Moreau theorem and the minimax theorem \cite[Corollary~37.3.2]{R70} (since $K$ is compact), that
\begin{equation}\label{Mimax}
\begin{aligned}
\max_{v \in K \cap \Im D^*_SN^*}\dotp{-(ND_S)^\dag ND_Te-\proj_Sv_0}{v}
&=\max_{v\in K}\pa{\dotp{-(ND_S)^\dag ND_Te-\proj_Sv_0}{v}-\iota_{\Im D^*_SN^*}(v)} \\
&=\max_{v\in K}\inf_{z \in \Ker N D_S}\dotp{-(ND_S)^\dag ND_Te-\proj_Sv_0+z}{v} \\
&=\min_{z \in \Ker N D_S}\max_{v\in K}\dotp{-(ND_S)^\dag ND_Te-\proj_Sv_0+z}{v} \\
&=\min_{z \in \Ker N D_S}\sigma_{K}(-(ND_S)^\dag ND_Te-\proj_Sv_0+z) .
\end{aligned}
\end{equation}
On the other hand, for any $v \in K \cap \Im D^*_SN^* \subset C^\circ$, we have $\sigma_C(v) \leq 1$ (see \eqref{Polar}). Let $v^\star \in K \cap \Im D^*_SN^*$ be a maximizer of the left hand side of \eqref{Mimax}. We then have from \eqref{gau} and \eqref{Mimax} that
\begin{equation}\label{eq:SIineq}
\begin{aligned}
0 \leq \min_{z \in \Ker N D_S}\sigma_{K}(-(ND_S)^\dag ND_Te-\proj_Sv_0+z)
&=-\dotp{(ND_S)^\dag ND_Te+ \proj_Sv_0}{v^\star} \\
&\leq \pa{1-\frac{c}{r\|D^*_S\|}}\sigma_C(v^\star) \leq 1-\frac{c}{r\|D^*_S\|} < 1 .
\end{aligned}
\end{equation}
It remains now to compute $\sigma_K$. We have from Fenchel-Moreau theorem, \cite[Corollary~16.4.1]{R70}, standard conjugacy calculus and \eqref{PoGa}, that
\begin{align*}
\sigma_K 
= \pa{\iota_{C^\circ} + \iota_{S}}^* 
= \cl\pa{\inf_{\eta \in \R^p} \sigma_{C^\circ}(\cdot+\eta) + \sigma_{S}(\eta)} 
= \cl\pa{\inf_{\eta \in T} \sigma_{C^\circ}(\cdot+\eta)} 
= \cl\pa{\inf_{\eta \in T} \gamma_{C}(\cdot+\eta)} .
\end{align*}
Since $\dom \gg_C = S$ by \cite[Propopsition~2]{VGFP15}, we have for any $z \in \R^p$ that
\[
\inf_{\eta \in T} \gamma_{C}(z + \eta) = \inf_{\eta \in T} \gamma_{C}(\proj_S z +\eta) = \gamma_{C}(\proj_S z) .
\]
Thus, the closure operation can be omitted above to get for any $z \in S$
\[
\sigma_K(z) = \gamma_C(z) .
\]
Inserting this into \eqref{eq:SIineq}, we get \ref{mtheo:claim2}, after observing that $\Im (ND_S)^\dag = \Im D_S^* N^* \subset S$ and that the minimum in $\rho(e)$ is achieved on $S$ (see Proposition~\ref{Sol}).

~[\ref{mtheo:claim2}$\Rightarrow$\ref{mtheo:claim3}]. Let $u^\star \in S$ be a solution of the minimization problem in \eqref{xyz}, denote $\ou = -(ND_S)^\dag NDe+u^* \in S$ and $\bar{v}=e+\ou$. We have $\bar{v}_T=e$ and, since \ref{mtheo:claim2} is satisfied, $\gg_{C}(\bar{v}_S-\proj_Sv_0)<1$. This is equivalent to $\bar{v} \in  \ri \partial J(x_0)$ thanks to the last claim of Lemma~\ref{Decom}. It remains to show that $D\bar{v} \in \Im \Phi^*$. As the Restricted Injectivity condition \eqref{GIR} holds, 
\begin{equation*}
ND_S(ND_S)^\dag ND_Te= ND_S (ND_S)^*(ND_S(ND_S)^*)^{-1}ND_Te=ND_Te,
\end{equation*}
which verifies \eqref{yz}. \eqref{xyz} is also verified as shown above and we can then argue as in the proof of the sufficient part in Proposition~\ref{Sol} to deduce that $D\bar{v} \in \Im \Phi^*$.

~[\ref{mtheo:claim3}$\Rightarrow$\ref{mtheo:claim1}]. Suppose that the Restricted Injectivity \eqref{GIR} holds at $x_0$ and there exists  $Dv\in \Im\Phi^*$ and $v \in \ri \partial J(x_0)$. Hence, $v_T=e$ and $\gg_{C}(\proj_S(v-v_0))<1$ thanks to Lemma~\ref{Decom}. Thus, for any $w\in \Ker \Phi$, we get from \eqref{dPhi} that 
\begin{eqnarray}\label{lower}\begin{array}{ll}
d\Psi(x_0)(w)
&=\dotp{e}{D^*w}+\dotp{\proj_Sv_0}{D^*_Sw}+\sigma_C(D^*_Sw)\\
&=\dotp{v-\proj_Sv}{D^*w}+\dotp{\proj_Sv_0}{D^*_Sw}+\sigma_C(D^*_Sw)\\
&=\dotp{-\proj_Sv}{ D^*w}+\dotp{\proj_Sv_0}{D^*_Sw}+\sigma_C(D^*_Sw)\\
&=\dotp{-\proj_S(v-v_0)}{ D^*_Sw}+\sigma_C(\proj_Sw)\\
&\ge -\sigma_{C}(D^*_Sw)\gg_{C}(\proj_S(v-v_0))+\sigma_C(D^*_Sw)\\
&=\Bpa{1-\gg_{C}(\proj_S(v-v_0))}\sigma_C(D^*_Sw) ,
\end{array}
\end{eqnarray}
where in the last two inequalities, we used the duality inequality on $\dom \sigma_C \times \dom \gamma_C=\R^p \times \dom S$. 
Recall that $0\in \ri C$ and $\B_s(0)\cap S\subset C$ for some $s > 0$. As $D^*_Sw\in S$, we have 
\[
\sigma_C(D^*_Sw)\ge \dotp{D^*_Sw}{s\frac{D^*_Sw}{\|D^*_Sw\|}}= s\|D^*_Sw\|\ge sc_1\|w\|,
\]
where $c_1 \eqdef \min\enscond{\|D^*_Sw\|}{w \in \mathbb{S}^{n-1} \cap \Ker\Phi}$, and $c_1 > 0$ by virtue of the Restricted Injectivity \eqref{GIR}. We derive from \eqref{lower} that
\[
d\Psi(x_0)(w) \ge \Bpa{1-\gg_{C}(\proj_S(v-v_0))}sc_1\|w\| \qforallq w\in \Ker \Phi,
\]
which verifies \ref{mtheo:claim1} by Lemma~\ref{Fa}. Recall that $Dv \in \Im \Phi^*$ and $v_T=e$, and thus $ND v_S = ND_T e$, which in turn shows that the vector $v_S$ obeys the constraint in \eqref{ro2}, and thus $\rho(e) \leq \gg_{C}(\proj_S(v-v_0)) < 1$. The sharpness constant of $\Psi$ can then be as large as devised in \eqref{ShaCo}. The proof is complete.
\end{proof} 
 
\begin{Corollary}[Sufficient condition for sharp solution] If the Restricted Injectivity holds at $x_0$ and the following condition
\begin{equation}\label{xi}
\tau(e) \eqdef \gg_C(-(ND_{S})^\dag NDe-\proj_Sv_0)<1
\end{equation}
is satisfied, then $x_0$ is a sharp solution to problem \eqref{BP0}.
\end{Corollary}
\begin{proof} 
It is easy to see that $\tau(e)\ge \rho(e)$. The result follows from Theorem~\ref{mtheo}. 
\end{proof} 

\begin{Remark}
The Restricted Injectivity \eqref{GIR} was proposed in \cite{FPVDS13,V14,VGFP15,VPF15}. It also has tracks in some special cases, e.g., for the $\ell_1$ norm problems \cite{F05,G17,GHS11,NDEG13,VPDF13,ZYC15, ZYY16}, $\ell_1/\ell_2$ norm problems \cite{G11,JKL15}, and nuclear norm problems \cite{CRPW12,CR09,CR13}. The form of the Nondegenerate Source Condition in Theorem~\ref{mtheo} appears also in \cite{V14,VGFP15,VPF15}. It generalizes the condition in \cite{GHS11} for the $\ell_1$ problem, which actually occurred earlier in \cite{CT04}. The combination of Restricted Injectivity and Nondegenerate Source Condition in the above result are proved in \cite[Theorem~5.3]{V14} as sufficient conditions for  solution uniqueness to problem \eqref{BP0}, but they are not necessary in general. By revisiting  Example~\ref{Ex12}, we see that $x_0$ is a unique solution to problem \eqref{Ex123} but the Nondegenerate Source Condition is not satisfied. Indeed, the only vector $v \in \R^2$ satisfying $\Phi^*v=(v_1+v_2,v_1,-v_2)^\top \in \partial J(x_0)=\{(0,1)\}\times[-1,1]$ is $v=(1,-1)$, but $\Phi^*v = (0,1,1) \notin \ri \partial J(x_0)$.

Theorem~\ref{mtheo} strengthens \cite[Corollary~1]{VGFP15} by showing Restricted Injectivity and Nondegenerate Source Condition are necessary and sufficient for sharp minima. In the case of $\ell_1$ problem, out result recovers  part of \cite[Theorem~2.1]{ZYY16} that gives  a characterization for solution uniqueness, which is equivalent to sharp minima in this framework. Some other characterizations are also studied recently in \cite{G17,MS19} by exploiting polyhedral structures. It is worth emphasizing that our result does not need polyhedrality. Our proof mainly  relies on the well-known first-order condition in Lemma~\ref{Fa} for sharp minima and the subdifferential decomposiability \eqref{Decop}. Theorem~\ref{mtheo} also covers many results about solution uniqueness in \cite{CRPW12, CR09,CR13, CT04, F05, G11, G17,GHS11, JKL15, NDEG13,RF08, VPDF13, T04, T06, ZYC15, ZYY16}.
\end{Remark}

%By Proposition~\ref{Sol}, $x_0$ is an optimal solution to problem \eqref{BP0}.  It suffices to show that $d \Psi(x_0;w)\ge \hat c\|w\|$ for any $w\in \Ker \Phi=\Im N^*$ with some constant $\hat c>0$ due to Fact~\ref{Fa}. Indeed, for any $u\in \R^r$ we get from \eqref{equ} and \eqref{mima} that
%\begin{eqnarray}\label{inq2}
%\begin{array}{ll}
%d\Psi(x_0;N^*u)&=\la ND_{S}(ND_{S})^\dag NDe,D^*_{S}N^*u\ra+\|D^*_{S}N^*u\|_\AA\\
%&\ge -\rho(e)\|D^*_{S}N^*u\|+\|D^*_{S}N^*u\|_\AA. 
%\end{array}
%\end{eqnarray}
%Since $\Ker D^*_{S}N^*=\{0\}$, there exists some constant $\al >0 $ such that $$\|D^*_{S}N^*u\|_\AA\ge \al\|u\|\ge\frac{\al}{\|N^*\|_{2,2}}\|N^*u\| \quad \mbox{for all}\quad u\in \R^r.$$
%This together with \eqref{inq2} tells us that 
%\[
%d\Psi(x_0;N^*u)\ge \frac{\al(1-\rho(e))}{\|N^*\|_{2,2}}\|N^*u\|\quad \mbox{ for all}\quad u\in \R^r,
%\]
%which verifies $d\Psi(x_0;w)\ge \hat c\|w\|$ for any $w\in \Ker \Phi$ with $\hat c \eqdef \frac{\al(1-\rho(e))}{\|N^*\|_{2,2}}>0$. The proof is completed. \end{proof} 

\begin{Remark}
The idea of using first-order analysis to study solution uniqueness to regularized optimization problems is not new as discussed above. For instance, the {\em Null Space Property} has been shown to ensure solution uniqueness for $\ell_1$ regularization \cite{DH01,NDEG13}. A generalization of this condition beyond the $\ell_1$ norm, coined {\em Strong Null Space Property}, was proposed in \cite{FPVDS13,V14,VGFP15}. This property reads 
\begin{equation}\label{SNS}
-\la e, D_T^*w\ra-\la \proj_Sv_0,D^*_Sw\ra<\sigma_C(D^*_{S}w) \qforallq w\in \Ker \Phi\setminus\{0\} .
\end{equation}
It is immediate to see from \eqref{Tes} that sharpness at $x_0$ entails \eqref{SNS}.
%However, our condition in (ii) is more quantitative and could be tested numerically; see Section~6. %Indeed, Proposition~\ref{Sol} suggests the way to compute $\rho(e)$ as the optimal value to
\end{Remark}

\paragraph{The case of analysis $\ell_1$}
When $\|\cdot\|_\AA=\|\cdot\|_1$, the model vector $e$ in \eqref{e} is indeed $(\sign(D_I^*x_0),0_K)$, where $I \eqdef \supp(D^*x_0) \eqdef \enscond{i\in\{1,2, \ldots,p\}}{(D^*x_0)_i \neq 0}$, $K \eqdef \{1, \ldots, p\}\setminus I$,  and $D_I$ is the submatrix of $D$ with column indices $I$. In his case, inequality \eqref{xi} takes the form
\[
{\tau}(\sign(D^*_Ix_0))=\|(ND_{K})^\dag ND_I\sign(D^*_Ix_0)\|_\infty<1,
\]
which is called the {\em Analysis Exact Recovery Condition} in \cite{NDEG13}. 
Another criterion used in \cite{VPDF13} to check solution uniqueness for $\ell_1$ problem is the so-called {\em Analysis Identifiability Criterion} at $\sign(D^*_Ix_0)$ denoted by
\begin{equation}\label{IC}
{\bf IC}(\sign(D^*_Ix_0)) \eqdef \min_{v\in \Ker  D_{J}} \|D_{K}^\dag(\Phi^*(\Phi^\dag_U)^*U^*-\Id)D_I \sign(D^*_Ix_0)-v\|_\infty<1,
\end{equation} 
where $U$ is basis of $\Ker D_{J}^*$ and $\Phi_U \eqdef \Phi U$. This condition reduces to the synthesis one introduced in \cite{F05} in the case of $D^*=\Id$ while ${\tau}(\sign(D^*_Ix_0))$ does not. As discussed in \cite{NDEG13,VPDF13}, $\tau(\sign(D^*_Ix_0))$ and ${\bf IC}(\sign(D^*_Ix_0))$ are different and no one implies the other even for the case $D=\Id$. 

For our general framework with $J$ as in \eqref{eq:Janalysis}, the Analysis Identifiability Criterion is satisfied at $x_0$ if 
\begin{equation}\label{AIC}
{\bf IC}(e) \eqdef \min_{u\in \Ker  D_S} \gg_C(D_S^\dag(\Phi^*(\Phi^\dag_U)^*U^*-\Id)D_Te-\proj_Sv_0+u)<1,
\end{equation} 
where $U$ is a matrix whose columns form a basis of $\Ker D^*_S$ and $\Phi_U \eqdef \Phi U.$

\begin{Proposition}\label{comp} 
If the Restricted Injectivity \eqref{GIR} holds at $x_0$ then $\rho(e) \le {\bf IC}(e)$. Consequently, if the Analysis Identifiability Criterion holds at $x_0$, then $x_0$ is a sharp solution to \eqref{BP0}.
\end{Proposition}
\begin{proof} 
\eqref{GIR} means that $\Ker \Phi_U=\{0\}$ and $\Phi_U^\dag=(\Phi_U^*\Phi_U)^{-1}\Phi_U^*$. Let $\ou\in\Ker D_S$ be a minimizer of problem \eqref{AIC}, and define   
\[
\ov \eqdef  D_S^\dag(\Phi^*(\Phi^\dag_U)^*U^*-\Id)D_Te+ \ou\quad \mbox{with}\quad \gg_C(\ov-\proj_Sv_0)={\bf IC}(e).
\] 
Note that 
 \[
 U^*(\Phi^*(\Phi^\dag_U)^*U^*-\Id)D_Te=(\Phi_U^*\Phi_U(\Phi_U^*\Phi_U)^{-1}U^*-U^*)D_Te=0.
 \]
It follows that $(\Phi^*(\Phi^\dag_U)^*U^*-\Id)D_Te\in \Im D_S,$ which implies that $D_S\ov=(\Phi^*(\Phi^\dag_U)^*U^*-\Id)D_Te$ by \eqref{cons}. Hence we have 
\begin{equation}\label{Nd}
ND_S\ov=(N\Phi^*(\Phi^\dag_U)^*U^*-N)D_Te=-ND_Te,
\end{equation}
which implies that $\ov\in -(ND_S)^\dag ND_Te+\Ker ND_S$ and thus $\rho(e) \le \gg_C(\ov-\proj_Sv_0)={\bf IC}(e)$. When the Analysis Identifiability Criterion is satisfied at $x_0$, we have $\rho(e)\le {\bf IC}(e)< 1$ and thus $x_0$ is a sharp solution to \eqref{BP0} due to Theorem~\ref{mtheo}. The proof is complete. 
\end{proof} 

In plain words, Proposition~\ref{comp} tells us that $\rho(e) < 1$ is weaker than the Analysis Identifiability Criterion \eqref{AIC}. In turn Theorem~\ref{mtheo} is stronger than \cite[Theorem~2]{VPDF13} for the analysis $\ell_1$ problem. It also covers the \cite[Proposition~5.7]{V14}. For the analysis $\ell_1$ problem, \cite{ZYY16} provides an example where $\rho(\sign(D^*_Ix_0))$ is  strictly smaller than both $\tau(\sign(D^*_Ix_0))$ and ${\bf IC}(\sign(D^*_Ix_0))$. %Moreover, by defining $\ov \eqdef (s_0, \ov_J)$, we get from \eqref{Nd} that $D\oy\in \Im \Phi^*$. Hence the Analysis Identifiability Criterion at $s_0$ means that $\oy\in {\rm ri}\,\partial \|D^*x_0\|_1 $, which means that 
%\begin{equation*}\label{NonD}
%D\ov\in {\rm ri}\, \partial J(x_0)\quad \mbox{with}\quad J(x)=\|D^*x\|_1,
%\end{equation*}
%which is known as the {\em non-degenerate} condition. It together with the Restricted Injectivity Condition has been used  in \cite{V14,VPF15} to prove robust recovery for more general convex regularizers $J(x)$ in \eqref{BP0}. However, this non-degenerate condition is nontrivial even for the case of $\ell_1$ norms.   

%%%%%%%%%%%%%%%%%%%%%%%%%%%
\subsection{Robust recovery with analysis decomposable priors}
The following result proves robust recovery with linear rate under Restricted Injectivity and Nondegenerate Source Condition.

\begin{Corollary}[Robust recovery of decomposable norm minimization]\label{Corode} Suppose that $J_0$ is a nonegative continuous convex function satisfying 
\begin{equation}\label{AC2}
    (D^*)^{-1}(\Ker (J_0)_\infty)\cap\Ker \Phi=\{0\}.
\end{equation}
If \eqref{GIR}-\eqref{NSC} hold at $x_0$, we have linear convergence rate for robust recovery as in \eqref{Est2} and \eqref{Est}. 
\end{Corollary}
\begin{proof} 
From \cite[Proposition~2.6.3]{AT03} $J_\infty(w)=(J_0)_\infty(D^*w)$ for any $w\in \R^n$. Hence, condition \eqref{AC2} is exactly \eqref{AC}. We obtain the claim by combining Theorem~\ref{mtheo} and Corollary~\ref{RobCon}.
\end{proof} 

The constant in $\mathcal{O}(\delta)$ can be made explicit. In particular the sharpness constant is given by \eqref{ShaCo}. This reveals that the "distance" to nondegeneracy is naturally captured by the sharpness constant, which plays a crucial role. Thus, the less degenerate, the more robust is the recovery. 
%If additionally, $\|\cdot\|_\AA$ is decomposable, the constant $s$ in \eqref{ShaCo} can be chosen as $1$ and $v_0=e_0$. This together with Theorem~\ref{Rob} allows us to obtain explicit bounds for robust recovery with linear rate.

%If $J(x)=\|D^*x\|_\AA$, by Theorem~\ref{mtheo}, $x_0$ is the unique solution to \eqref{BP0}. By the discussion in Remark~\ref{Relax}, condition \eqref{AC} is satisfied. The proof is complete. \end{proof} 

Corollary~\ref{Corode} covers parts of \cite[Theorem~2]{FPVDS13} (see also \cite{V14,VPF15}), \cite[Theorem~4.2]{GHS11}, and \cite[Theorem~2]{ZYY16} for the case of (analysis or synthesis) $\ell_1$ minimization problems. When $D^*=\Id$ and $J_0(x)=\|x\|_\AA$, it also covers the results in \cite{CRPW12} about robust recovery for the constrained problem \eqref{BP1}.

%One natural question arising from Corollary~\ref{ell1} is: 
%\begin{enumerate}[{\bf (Q1)}]
%\item Can a unique solution be sharp to a 
%decomposable norm regularized problem?
%\end{enumerate}
%Except for the case of $\ell_1$ norm as discussed in Corollary~\ref{ell1}, we have negative answers for both group sparsity and low-rank problems; see our Example~\ref{Last3}, Example~\ref{Ex12}, and Example~\ref{Last}. However,  the next section shows  that a unique solution to $\ell_1/\ell_2$ problem  is indeed the strong solution. %For nuclear norm problem, the answer for {\bf (Q1)} is negative too. Particularly,  Example~\ref{Last} shows that $x_0$ is a unique solution to \eqref{BP0}, but it is neither  strong nor  sharp. 

%%%%%%%%%%%%%%%%%%%%%%%%%%%
\subsection{Connections between unique/sharp/strong solutions in the noiseless case}
When there is noise in observation \eqref{yw}, problem \eqref{Las3} is usually used to recover the original signal $x_0$. Solution uniqueness to \eqref{Las3} is especially important for exact recovery \cite{F05,VPDF13}. We show next that Restricted Injectivity and Nondegenerate Source Condition are sufficient for strong minima to problem \eqref{Las3}, and they become necessary for the $\ell_1$ problem. This result is true for a larger class of regularized problems taking the form
\begin{equation}\label{Las0}
\min_{x\in \R^n}\quad \Theta(x) \eqdef f(\Phi x)+\mu J(x),
\end{equation}
where $\mu$ is a positive parameter and the loss function $f:\R^m\to \oR$ is an extended real-valued convex function satisfying the following two conditions:
\begin{enumerate}[label={\rm (\Alph*)}]
\item $f$ is twice continuously  differentiable in $\Int (\dom f)$. \label{assum:A}
\item $\nabla f^2(x)$ is positive definite for all $x\in \Int (\dom f)$, i.e., $f$ is strictly convex in the interior of its domain. \label{assum:B}
\end{enumerate}
In \eqref{Las3}, the function $f$ is $\frac{1}{2}\|\cdot-y\|^2$ which certainly satisfies the above two conditions. Moreover, the standing assumptions \ref{assum:A} and \ref{assum:B} for $f$ cover the important case of the Kullback-Leiber divergence:
\begin{eqnarray}\label{KL}
f_{\rm KL}(z)=\left\{\begin{array}{ll} \disp\sum_{i=1}^my_i\log\dfrac{y_i}{z_i}+z_i-y_i\quad &\mbox{if}\quad z\in \R^m_{++}\\
 +\infty \quad &\mbox{if}\quad z\in \R^m_+\setminus\R^m_{++},\end{array}\right.
\end{eqnarray}
where $y \in \R^m_{+}$ and $0 \log 0 = 0$. This offers a natural way to measure of similarity of two nonnegative vectors (e.g., two discrete distributions) and is broadly used in statistical/machine learning and signal processing.% or the loss information when using prior distribution to approximate the posterior.

In the following result, we provide the connections between unique/sharp/strong solutions for the two problems \eqref{BP0} and \eqref{Las0}. Part~\ref{Stheo:claim1} of this result could be obtained from \cite[Proposition~3.2]{MS19}, but we still give  a short proof as our assumptions on $f$ are slightly different, e.g., $f$ may not have full domain.  

\begin{Proposition}[Unique/sharp/strong solutions for problems \eqref{BP0} and \eqref{Las0}]\label{Stheo} 
Suppose that $\ox$ is an optimal solution to problem \eqref{Las0} and $\Phi\ox\in \Int (\dom f)$. Then $\ox$ is a solution to problem \eqref{BP0} with $x_0=\ox$. Moreover, the following statements hold:
\begin{enumerate}[label={\rm (\roman*)}]
\item $\ox$ is the unique solution to \eqref{Las0} if and only if it is the unique solution to \eqref{BP0}. \label{Stheo:claim1} 
\item If $\ox$ is a sharp solution to \eqref{BP0} then it is a strong solution to problem \eqref{Las0}. \label{Stheo:claim2} %Thus conditions $\Ker \Phi\cap \Ker D_{S}^*=\{0\}$ and  $\rho(e)<1$ are sufficient for strong solution at $\ox$ to problem \eqref{Las0}.
\item If $\ox$ is a strong solution to \eqref{Las0} then it is also a strong solution to \eqref{BP0}. \label{Stheo:claim3} 
\end{enumerate}
\end{Proposition}
\begin{proof} 
For any $x\in \Phi^{-1}(\Phi\ox)$, we have 
\begin{equation}\label{indy}
f(\Phi \ox)+\mu J(x)=f(\Phi x)+\mu J(x) \ge f(\Phi \ox)+\mu J(\ox).
\end{equation}
It follows that $\ox$ is a solution to problem \eqref{BP0} with $y_0=\Phi \ox$. 

\begin{enumerate}[label={\rm (\roman*)}]
\item If $\ox$ is the unique solution to problem \eqref{Las0}, we have the strict inequality in \eqref{indy} provided that $x\neq \ox$, which shows that $\ox$ is also the unique solution to problem \eqref{BP0}. 
%To get the converse claim, it is sufficient to show that $\ox$ is an isolated point to the convex solution set of problem \eqref{Las0} provided that  $\ox$ is the unique solution to \eqref{BP0}. 
Assume now that $\ox$ is the unique solution to \eqref{BP0}. By contradiction, suppose that $\ox$ is not the unique solution to \eqref{Las0}. Since $\Argmin(\Theta)$ is convex and $\Phi^{-1}(\Int(\dom f))$ is an open set containing $\ox$, there exists $r > 0$ such that $\B_r(\bar x) \subset \Phi^{-1}(\Int(\dom f))$ and $\B_r(\bar x) \cap \Argmin(\Theta) \neq \emptyset$, while the later is not the singleton $\{\ox\}$. Choose $\hat x \in \B_r(\bar x)$ with $\hat x \neq \ox$. We have $-\Phi^*\nabla f(\Phi\hat x)\in \partial J(\hat x)$ and $-\Phi^*\nabla f(\Phi\ox)\in \partial J(\ox)$. The monotonicity of the subdifferential tells us that 
\[
-\dotp{\nabla f(\Phi \hat x)-\nabla f(\Phi \ox)}{\Phi\hat x-\Phi\ox}=-\dotp{\Phi^*\nabla f(\Phi \hat x)-\Phi^*\nabla f(\Phi \ox)}{\hat x-\ox} \ge 0.
\]
Convexity of $f$ entails the opposite inequality which shows that
\[
\dotp{\nabla f(\Phi \hat x)-\nabla f(\Phi \ox)}{\Phi\hat x-\Phi\ox}=0.
\]
Since $f$ is twice continuously differentiable in $\Int (\dom f)$, we obtain from the mean-value theorem that 
\begin{equation}\label{nabla}
\int_0^1\dotp{\nabla^2 f(\Phi(\ox+t(\hat x-\ox)))\Phi(\hat x-\ox)}{\Phi(\hat x-\ox)} dt = 0.
\end{equation}
Since the Hessian is positive definite and continuous on $\Int (\dom f)$, there exists some $\tau>0$ such that 
\[
\dotp{\nabla^2 f(\Phi(\ox+t(\hat x-\ox)))\Phi(\hat x-\ox)}{\Phi(\hat x-\ox)} \ge \tau \|\Phi(\hat x-\ox)\|^2
\]
for all $t\in [0,1]$.
Combining this with \eqref{nabla} implies that $\Phi\hat x=\Phi\ox$. Using this together with the fact that both $\hat x$ and $\ox$ are minimizers to \eqref{Las0}, entails that $J(\hat x)=J(\ox)$. This contradicts uniqueness of $\ox$ for \eqref{BP0}.

\item Assume that $\ox$ is a sharp solution to \eqref{BP0}. By Proposition~\ref{Shar}, we have 
\begin{equation}\label{DD}
\Ker \Phi\cap \CC_J(x_0)=\{0\}
\end{equation}
where we recall $\CC_J(x_0)$ from \eqref{D0}. From the sum rule \eqref{sum} and convexity of $J$, we have 
\begin{equation}\label{Tep}
d^2\Theta(\ox|0)(w)=\dotp{\nabla^2 f(\Phi \ox)\Phi w}{\Phi w}+\mu d^2J (\ox|-\mu^{-1}\Phi^*\nabla f(\Phi\ox)(w))\ge \dotp{\nabla^2 f(\Phi \ox)\Phi w}{\Phi w}
\end{equation}
and from \eqref{dom2}, $\dom d^2\Theta(\ox|0)\subset \enscond{w\in \R^n}{d\Theta(\ox)(w)= 0}$. We also have 
\[
d\Theta(\ox)(w)=\dotp{\nabla f(\Phi\ox)}{\Phi w}+\mu dJ(x_0)(w).
\]
It then follows that 
\begin{equation}\label{Kx}
\dom d^2\Theta(\ox|0)\subset\enscond{w\in \R^n}{\dotp{\nabla f(\Phi\ox)}{\Phi w}+\mu dJ(x_0)(w)=0}.
\end{equation}
To verify that $\ox$ is a strong solution to problem \eqref{Las0} by using Lemma~\ref{Fa}, we claim that  
\begin{equation}\label{Ky}
\dotp{\nabla^2f(\Phi \ox)\Phi w}{\Phi w} >0 \qforallq w \in \dom d^2\Theta(\ox|0)\setminus\{0\}.
\end{equation}
Suppose that  $w\in \dom d^2\Theta(\ox|0)$ satisfying $\dotp{\nabla^2f(\Phi \ox)\Phi w}{\Phi w} \le 0$. Since   $\nabla^2f(\Phi \ox) \succ 0$, we have $\Phi w=0$. This together with \eqref{Kx} tells us that $w \in \Ker \Phi\cap \CC_J(x_0)=\{0\}$ by \eqref{DD}. Thus inequality \eqref{Ky} holds and $\ox$ is a strong solution to problem  \eqref{Las0}. 

\item Suppose that $\ox$ is a strong solution to problem \eqref{Las0}. There exist constants $\kk,\gg>0$ such that
\[
f(\Phi x)+\mu J(x)\ge f(\Phi\ox)+\mu J(\ox)+\frac{\kk
}{2}\|x-\ox\|^2 \qforallq x\in \B_\gg(\ox).
\]
For any $x\in \B_\gg(\ox)\cap \Phi^{-1}(\Phi \ox)$, we obtain that 
\[
 J(x)\ge J(\ox)+\frac{\kk}{2\mu}\|x-\ox\|^2,
\]
which means that $\ox$ is also a strong solution to \eqref{BP0} as claimed. 
\end{enumerate}
\end{proof}

For the special case \eqref{Las3} of \eqref{Las0}, \cite[Theorem~1]{FPVDS13} shows that if $\ox$ is an optimal solution to \eqref{Las0} with $x_0=\ox$  and the Strong Null Space Property \eqref{SNS} holds at $\ox$, $\ox$ is the unique  solution to \eqref{Las0}. As Strong Null Space Property is a characterization for sharp minima to problem \eqref{BP0},  Proposition~\ref{Stheo} advances \cite[Theorem~1]{FPVDS13} and \cite[Theorem~3]{VGFP15} with further information that $\ox$ is a strong solution to \eqref{Las3}.

 %This result could be explained as follows. Since the Strong Null Space Property is a characterization for sharp minima at $\ox=x_0$ to problem \eqref{BP0} as discussed before \eqref{SNS}, $\ox$ is the unique solution to problem \eqref{BP0}. If $\hat x$ is another solution to \eqref{Las0}, we have $-\Phi^*(\Phi \hat x-y)\in D\partial\|D\hat x\|_\AA$ and $-\Phi^*(\Phi \ox-y)\in D\partial\|D\hat x\|_\AA$, due to the monotonicity of the subdifferential, one has
 %\[
 %\la -\Phi^*(\Phi \hat x-y)+\Phi^*(\Phi \ox-y),\hat x-\ox\ra\ge 0,
% \]
% which implies $\Phi\hat x=\Phi \ox$. Thus $\hat x$ is also a solution to problem \eqref{BP0}, which means $\hat x=\ox$, i.e., $\ox$ is the unique solution to \eqref{Las0}. Another direct explanation could come from \cite[Proposition~3.2]{MS19}. In the following result, we provide another insight, which tells us that Strong Null Space Property or our conditions are sufficient for strong solution to problem \eqref{Las0}. 

Two natural questions arise from Proposition~\ref{Stheo}: are the converse statements of \ref{Stheo:claim2}-\ref{Stheo:claim3} true ? That is:
\begin{enumerate}[label={\rm \bf (Q.\arabic*)}]
\item \label{OQ1}
If $\ox$ is a strong solution to \eqref{Las0}, can it be a sharp solution to \eqref{BP0} with $x_0=\ox$? 
\item \label{OQ2}
If $\ox$ is a strong solution to \eqref{BP0}, can it be a strong solution to \eqref{Las0}? 
\end{enumerate}
For analysis $\ell_1$ problems, i.e., $J_0(\cdot)=\|\cdot\|_1$, and more generally for $J_0$ the support function of any polyhedral convex compact $C$ such that $0 \in \ri C$,  we have positive answers for both questions.
\begin{Corollary}[Solution uniqueness to $\ell_1$ problems]\label{ell2} 
Let $\ox$ be a minimizer to \eqref{Las0} with $J_0(\cdot)=\|\cdot\|_1$. Then the following are equivalent:
\begin{enumerate}[label={\rm (\roman*)}]
\item $\ox$ is the  unique solution to problem \eqref{BP0} with $x_0=\ox$. \label{ell2:claim1}
\item $\ox$ is the sharp solution to problem \eqref{BP0} with $x_0=\ox$. \label{ell2:claim2}
\item $\ox$ is the unique solution to problem \eqref{Las0}. \label{ell2:claim3}
\item $\ox$ is the  strong solution to problem  \eqref{Las0}. \label{ell2:claim4}
\end{enumerate}
\end{Corollary}
\begin{proof} 
~[\ref{ell2:claim1}$\Leftrightarrow$\ref{ell2:claim3}] and [\ref{ell2:claim2}$\Rightarrow$\ref{ell2:claim4}] are from Proposition~\ref{Stheo}. [\ref{ell2:claim4}$\Rightarrow$\ref{ell2:claim3}] is trivial. Finally, [\ref{ell2:claim1}$\Leftrightarrow$\ref{ell2:claim2}] follows by combining Proposition~\ref{DeUn} and Proposition~\ref{Shar}\ref{Shar:claim1} since $J_0$ is a polyhedral norm, in which case the descent and critical cones coincide at any $x_0$ (see Proposition~\ref{prop:tancritcones}). 
 \end{proof}

According to Theorem~\ref{mtheo}, conditions $\Ker \Phi\cap \Ker D_{S}^*=\{0\}$ 
and  $\rho(e)<1$ form a characterization for solution uniqueness to $\ell_1$ problem.  A similar result was established in \cite[Theorem~1]{ZYY16} for problems \eqref{BP0} and \eqref{Las3} with the extra assumption that $\Phi$ has full row-rank. Corollary~\ref{ell2} is more general and reveals that the unique solution to problem \eqref{Las0} is the strong solution.  

The answer is negative for \ref{OQ1} in general; see our Theorem~\ref{Las2} and Example~\ref{Ex12}. Regarding \ref{OQ2}, we do have a positive answer for group-sparsity ; see Theorem~\ref{BP12} and Theorem~\ref{Las2}. However, it is not true for the nuclear norm miniminzation problem as we now show.   

\begin{Example}[Strong solutions of \eqref{BP0} are not those of \eqref{Las0}]\label{Last3}
Let us consider the following nuclear norm minimization problem
\begin{equation}\label{LN}
    \min_{X\in\R^{2\times 2}} \Theta(X)  \eqdef  \dfrac{1}{2}\|Y - \Phi(X)\|^2+ \|X\|_*,
\end{equation}
where $\|\cdot\|_*$ stands for the nuclear norm of $X$, and $\Phi:\R^{2\times2}\to \R^{2}$ is the diagonal operator, i.e., $\Phi(X)=(X_{11},X_{22})^\top$ and $Y=(2,1)^\top$. This is a special case of \eqref{Las0} with $f(\cdot)=\frac{1}{2}\|Y-\cdot\|^2$. Let 
$\Bar X=
\begin{pmatrix}
1   &  0\\
0  &  0
\end{pmatrix}$. 
We have 
\begin{equation*}
\Phi^*\nabla f(\Phi(\Bar X)) = 
\begin{pmatrix}
  -1   &  0\\
   0  &  -1
\end{pmatrix} \qandq \partial\|\Bar X\|_* = 
\enscond{\begin{pmatrix}
  1   &  0\\
  0  &  \alpha
\end{pmatrix}}{\alpha\in [-1,1]}.
\end{equation*}
Thus $0\in \partial \Theta(\Bar X)$, i.e., $\Bar X$ is a solution to \eqref{LN}. By Proposition~\ref{Stheo}, it is also a solution to
\begin{equation}\label{BPN}
\min_{X\in\R^{2\times 2}} \|X\|_* \qandq \Phi(X)=(1,0)^\top.
\end{equation}
For any $2\times 2$ matrix $X$, let $\sigma_1,\sigma_2$ be its singular values. We have
\begin{equation}\label{nucl}
   \|X\|_*=\sigma_1+\sigma_2=\sqrt{\sigma_1^2+\sigma^2_2+2\sigma_1\sigma_2}=\sqrt{\|X\|_F^2+2|{\rm det}\,(X)|}. 
\end{equation}
The feasible set of \eqref{BPN} consists of matrices of the form 
$X= \begin{pmatrix}
1   &  a\\
b  &  0
\end{pmatrix}$, and we obtain from \eqref{nucl} that 
\[
\|X\|_*= \sqrt{1+(|a|+|b|)^2}.
\]
Thus $\Bar X$ is the unique solution to \eqref{BPN} as well as $\eqref{LN}$ due to Proposition~\ref{Stheo} again. Next we claim that $\Bar X$ is the strong solution to \eqref{BPN}, which means there exist $\kappa, \delta>0$ such that
\begin{equation}\label{nucl_str}
\|X\|_*-\|\Bar X\|_*\geq\dfrac{\kappa}{2}\|X-\Bar X\|^2_F \qforallq   X= \begin{pmatrix}
1   &  a\\
b  &  0
\end{pmatrix} \qwithq \|X-\Bar X\|_F=\sqrt{a^2+b^2} \le \delta.
\end{equation} 
Indeed, we have 
\begin{align*}
\|X\|_*-\|\Bar X\|_* 
&= \sqrt{1+(|a|+|b|)^2}-1\\
&=\dfrac{a^2+b^2+2|ab|}{\sqrt{1+a^2+b^2+2|ab|}+1}\\
&\geq \dfrac{a^2+b^2}{\sqrt{1+2\delta^2}+1}\\
&= \dfrac{1}{\sqrt{1+2\delta^2}+1}\|X-\bar X\|_F^2.
\end{align*}
This certainly verifies \eqref{nucl_str}. Next we claim that $\Bar X$ is not a strong solution to \eqref{LN}. Pick $X_\ve \eqdef \begin{pmatrix}
1+\ve^2 & \ve \\
 \ve    & \ve^2
\end{pmatrix}$ with $\ve>0$ sufficiently small, observe that $\|X-\Bar X\|^2_F=2(\ve^2+\ve^4)$. It follows from \eqref{nucl} that
\begin{align*}
\Theta(X_\ve)-\Theta(\Bar X) 
&= \dfrac{1}{2}
\anorm{\begin{pmatrix}
    \ve^2-1   \\
    \ve^2-1
\end{pmatrix}}^2+
\anorm{\begin{pmatrix}
    1+\ve^2 & \ve   \\
    \ve & \ve^2
\end{pmatrix}}_*-2\\
&= (\ve^2-1)^2+\sqrt{((1+\ve^2)^2+2\ve^2+\ve^4)+2\ve^4}-2\\
&= (\ve^2-1)^2+(1+2\ve^2)-2\\
&=\ve^4.
\end{align*}
Therefore, $\Bar X$ cannot be a strong solution to \eqref{LN}, which also implies that $\Bar X$ is neither a sharp solution to \eqref{LN} nor to \eqref{BPN} due to Proposition~\ref{Stheo}.
%\color{red}{My calculation:
%\begin{eqnarray}
%d^2\Theta(\Bar X|0)(W)=\left\{\begin{array}{ll}
%a^2+d^2+(b-c)^2\quad &\mbox{if}\quad d>0,\\
%a^2+(|b|+|c|)^2\color{blue}{\longleftarrow  a^2+(b-c)^2????}\quad &\mbox{if}\quad d=0\\
%\infty \quad &\mbox{if}\quad d<0\end{array}\right.\quad \mbox{with}\quad W=\begin{pmatrix}
%a & b \\
%c    & d
%\end{pmatrix}.
%\end{eqnarray}
%It follows that $d^2\Theta(\Bar X|0)(W)>0$ for all $W\neq 0$, which tells us that $\Bar X$ is a strong solution!!!
%}
\end{Example}

%% file: tex/sec_strong_group.tex
\section{Characterizations of unique/strong solutions for group-sparsity}
\label{sec:stronggroup}
%%%%%%%%%%%%%%%%%%%%%%%%%%%%%%%%%%%%%%%%%%%%%%%%%%%%%%%%
In this section, we study the following particular case of \eqref{BP0}, where $J_0$ is the $\ell_1/\ell_2$ norm that promotes group sparsity \cite{RRN12,RF08,YL06}. 
\begin{equation}\label{ITV}
\min_{x\in \R^n}\quad \|D^*x\|_{\ell_1/\ell_2}\quad\mbox{subject to}\quad \Phi x=\Phi x_0 .
\end{equation}

Following the notation in \cite{CR13}, we suppose that $\R^p$ is decomposed into $q$ groups by 
\begin{equation}
    \R^p=\bigoplus_{g=1}^q V_g,
\end{equation}
where each $V_g$ is a subspace of $\R^p$ with the same dimension $G$. For any $u\in \R^p$, we write  $u=\disp\sum_{g=1}^qu_g$ with $u_g\in V_g$ being  the vector in group $V_g$. The $\ell_1/\ell_2$ norm in $\R^p$  is defined by
\begin{equation}\label{l12}
\|u\|_{\ell_1/\ell_2}=\sum_{g=1}^q\|u_g\|.  
\end{equation}
Its dual is the  $\ell_\infty/\ell_2$ norm:
\begin{equation}\label{li2}
\|u\|_{\ell_\infty/\ell_2}=\max_{1\le g\le q}\|u_g\|.  
\end{equation}
With  $\ou \eqdef D^*x_0$, define $I \eqdef \left\{g\in \{1,\ldots,q\}|\; \ou_g\neq 0\right\}$, the index set of active groups of $\ou$, and $K \eqdef \{1,\ldots,q\}\setminus I$, the index set of nonactive groups of $\ou$. Note that 
\begin{equation}\label{sub12}
    \partial\|\ou\|_{\ell_1/\ell_2} \eqdef \left\{v\in \R^p|\; v_g= \frac{ \ou_{g}}{\|\ou_{g}\|_2}\mbox{ for } g\in I \mbox{ and } \|v_g\|\le 1 \mbox{ for } g\in K\right\}. 
\end{equation}
Thus $\ell_1/\ell_2$ norm is decomposable at $\ou$ as in Remark~\ref{Deco}, where $T=\disp\bigoplus_{g\in I} V_g$, $S=T^\perp$, and 
\begin{equation}\label{e12}
    e=\sum_{g\in I} \frac{\ou_{g}}{\|\ou_{g}\|}.
\end{equation}

%%%%%%%%%%%%%%%%%%%%%%%%%%%
\subsection{Descent cone of group sparsity}
Sharp minima at $x_0$ for problem \eqref{ITV} is studied in our previous section. However, unlike the case of $\ell_1$ problem, a unique solution to problem \eqref{ITV} may be not sharp; see  Example~\ref{Ex12}. To characterize the solution uniqueness to group sparsity problem \eqref{ITV}, we compute the  descent cone $\DD_J(x_0)$ in \eqref{Des} as follows.  

\begin{Theorem}[Descent cone to $\ell_1/\ell_2$ problem and geometric characterization for solution uniqueness]\label{Uniq} The descent cone  at $x_0$ to problem \eqref{ITV} is given by 
\begin{equation}\label{Sub}
\DD_{J} (x_0)=(\EE\cap \bd \CC_{J}(x_0))\cup(\Int \CC_{J}(x_0)),
\end{equation}
where 
\begin{equation}\label{C} 
\EE \eqdef \enscond{w\in\R^n}{D_T^*w\in \Im\{(D^*_Tx_0)_g|\; g\in I\}} ,
\end{equation}
and $\bd \CC_{J}(x_0)$ stands for the boundary of the critical cone $\CC_{J}(x_0)$ defined in \eqref{D0}, i.e., 
\begin{equation}\label{bd}
\bd \CC_{J}(x_0)=\enscond{w\in \R^p|}{\dotp{D_Te}{w}+\|D^*_{S}w\|_{\ell_1/\ell_2}=0}.
\end{equation}
Consequently, $x_0$ is a unique solution to \eqref{ITV} if and only if 
\begin{equation}\label{Geo}
\Ker \Phi \cap \Bpa{(\EE\cap \bd \CC_{J}(x_0))\cup(\Int \CC_{J}(x_0))}=\{0\}.
\end{equation}
\end{Theorem}
\begin{proof} 
Let us start by verifying the inclusion ``$\supset$'' in \eqref{Sub}.  Recall that  $J(x)=\|D^*x\|_{\ell_1/\ell_2}$.  For any $w\in \Int (\CC_{J}(x_0))$, we get from \eqref{D0} that  
\[
dJ(x_0)(w)=\lim_{t\dn 0}\dfrac{J(x_0+tw)-J(x_0)}{t}<0.
\]
Hence there is some $t_0>0$ such that $J(x_0+t_0w)< J(x_0)$. This ensures $w\in \DD_J(x_0)$. For any $w\in \EE\cap \bd \CC_{J}(x_0)$, we represent $D_T^*w=\sum_{g\in I}\lm_g (D^*x_0)_g$ with $\lm_g\in \R$, $g\in I$ and 
% \begin{equation}\label{inqe}
%     \la D_Te,w\ra+\|P_{S}w\|_{\ell_1/\ell_2}=\la e,D^*_Tw\ra+\|P_{S}w\|_{\ell_1/\ell_2}=\sum_{g\in I}{\lm_g}+\sum_{g\in J}\|w_g\|= 0.
% \end{equation}
% Define $t_1 \eqdef \min\left\{-\dfrac{\|(D^*x_0)_g\|}{\lm_g}|\; g\in I_-\right\}>0$ with $I_- \eqdef \{g\in I|\; \lm_g<0\}$, we have 
% \[\begin{array}{ll}
% \|D^*x_0+t_1D^*w\|_{\ell_1/\ell_2}-\|D^*x_0\|_{\ell_1/\ell_2}&=\disp\sum_{g\in I}\left|\|(D^*x_0)_g\|+t_1\lm_g\right|+\sum_{g\in J} t_1\|D^*w_g\|-\|(D^*x_0)_g\|\\
% &=\disp\sum_{g\in I}\|(D^*x_0)_g\|+t_1\lm_g+\sum_{g\in J} t_1\|(D^*w)_g\|-\|(D^*x_0)_g\|,
% \end{array}
% \]
\begin{equation}\label{inqep}
 0 = \dotp{ D_Te}{w}+\|D^*_{S}w\|_{\ell_1/\ell_2}=\dotp{e}{D^*_Tw}+\|D^*_{S}w\|_{\ell_1/\ell_2}=\sum_{g\in I}{\lm_g}\|(D^*x_0)_g\|+\|D_{S}^*w\|_{\ell_1/\ell_2}.
\end{equation}
Define  $I_- \eqdef \{g\in I|\; \lm_g<0\}$. If $I_-=\emptyset$, we get from \eqref{inqep} that $D^*_Sw=0$ and $\lm_g=0$ for all $g\in I$, which implies  that $D^*_Tw=0$ and thus $D^*w=0$. Therefore, we get
\[
J(x_0+w)=J_0(D^*(x_0+w))=J_0(D^*x_0)=J(x_0),
\]
which tells us that $w\in \DD_J(x_0)$. If $I_-\neq\emptyset$, define
$t_1 \eqdef \min\enscond{-\dfrac{1}{\lm_g}|}{g\in I_-}>0$, we obtain from \eqref{inqep} that  
\begin{align*}
J(x_0+t_1w)-J(x_0)
&=\|D^*x_0+t_1D^*w\|_{\ell_1/\ell_2}-\|D^*x_0\|_{\ell_1/\ell_2} \\
&=\sum_{g\in I}|1+t_1\lm_g|\|(D^*x_0)_g\|+t_1\|D_{S}^*w\|_{\ell_1/\ell_2}-\|D^*x_0\|_{\ell_1/\ell_2}\\
&=\sum_{g\in I}(1+t_1\lm_g)\|(D^*x_0)_g\|+t_1\|D_{S}^*w\|_{\ell_1/\ell_2}-\|D^*x_0\|_{\ell_1/\ell_2}\\
&=t_1\pa{\sum_{g\in I}\lm_g\|(D^*x_0)_g\|+\|D_{S}^*w\|_{\ell_1/\ell_2}}+\sum_{g\in I}\lm_g\|(D^*x_0)_g\|-\|D^*x_0\|_{\ell_1/\ell_2}\\
&= 0,
\end{align*}
where the third equality is from the choice of $t_1$. In both cases of $I_-$,  $w\in \DD_J(x_0)$.

To justify the reverse inclusion ``$\subset$'' in \eqref{Sub}, take any $w\in \DD_J(x_0)$ and let $t>0$ such that $\|D^*x_0+tD^*w\|_{\ell_1/\ell_2}\le \|D^*x_0\|_{\ell_1/\ell_2}$. For any $\al\in (0,1)$, we have
\begin{equation}\label{Dxl}
\|D^*x_0\|_{\ell_1/\ell_2}\ge \al \|D^*x_0+tD^*w\|_{\ell_1/\ell_2}+(1-\al)\|D^*x_0\|_{\ell_1/\ell_2}\ge \|D^*x_0+\al tD^*w\|_{\ell_1/\ell_2}.
\end{equation}
Choose $\al>0$ sufficiently small such that 
\[
\|(D^*x_0)_g+\al t(D^*w)_g\|>0 \qforallq  g\in I. 
\]
Since $w\in \DD_J(x_0) \subset \CC_{J}(x_0)$ by \eqref{CDinc}, it suffices to show that if  $w\in \bd \CC_{J}(x_0)$ then $w\in \EE.$ Suppose that $w\in \bd \CC_{J}(x_0)$, we have 
\[
0=\dotp{D_Te}{w}+\|D_{S}^*w\|_{\ell_1/\ell_2}=\dotp{e}{D^*_Tw}+\|D_{S}^*w\|_{\ell_1/\ell_2}.
\]
It follows from the latter and \eqref{Dxl} that 
\begin{equation}\label{Big}
\begin{aligned}
0
&\ge\|D^*x_0+\al tD^*w\|_{\ell_1/\ell_2}- \|D^*x_0\|_{\ell_1/\ell_2}\\
&=\sum_{g\in I}\|(D^*x_0)_g+\al t (D^*w)_g\|+\al t\sum_{g\in K}\|(D^*w)_g\|-\|D^*x_0\|_{\ell_1/\ell_2}\\
&=\sum_{g\in I}\Bpa{\|(D^*x_0)_g+\al t (D^*w)_g\|-\|(D^*x_0)_g\|-\al t\dotp{e_g}{(D^*w)_g}}.
\end{aligned}
\end{equation}
Since $e_g= \partial \|(D^*x_0)_g\|$, each $\|(D^*x_0)_g+\al t (D^*w)_g\|-\|(D^*x_0)\|_g-\al t\dotp{e_g}{(D^*w)_g}\ge 0$ for $g\in I$. This together with \eqref{Big} tells us that
\[
0=\|(D^*x_0)_g+\al t (D^*w)_g\|-\|(D^*x_0)_g\|-\al t\dotp{e_g}{(D^*w)_g} \qforallq g\in I.
\]
Hence, we have
\[
0=\|(D^*x_0)_g\|-\dotp{e_g}{(D^*x_0)_g}=\|(D^*x_0)_g+\al t (D^*w)_g\|-\dotp{e_g}{(D^*x_0)_g+\al t (D^*w)_g}.
\]
As $\|e_g\|=1$ and $(D^*x_0)_g+\al t (D^*w)_g\neq 0$, the latter equality holds when
\[
(D^*x_0)_g+\al t (D^*w)_g=\delta_g e_g 
\]
for some $\delta_g>0$, $g\in I$. It follows that 
\[
(D^*w)_g=\frac{1}{\al t}(\delta_g -\|(D^*x_0)_g\|)e_g
\]
for any $g\in I$, which ensures that $w\in \EE$ and verifies the equality \eqref{Sub}. 

The characterization for solution uniqueness at $x_0$ in \eqref{Geo} follows directly from Proposition~\ref{DeUn} and \eqref{Sub}. \end{proof}

%%%%%%%%%%%%%%%%%%%%%%%%%%%
\subsection{Unique vs strong solutions}
We show next that a unique solution to \eqref{BP0} is indeed a strong solution. The proof is based on the second-order analysis in Lemma~\ref{Fa}. We need the computation of the second subderivative for the function $\Psi$ defined in \eqref{Psi}.

\begin{Lemma}[Second subderivative to $\ell_1/\ell_2$ norm]\label{SOD}  Suppose that $x_0$ is a minimizer to problem \eqref{ITV}. Then we have 
\begin{equation}\label{dom}
    \dom d^2\Psi(x_0|0)=\Ker \Phi\cap \bd \CC_{J}(x_0)
\end{equation} 
and  
\begin{eqnarray}\label{d2}
d^2\Psi(x_0|0)(w)=\sum_{g\in I} \dfrac{\|(D^*w)_g\|^2\|(D^*x_0)_g\|^2-\la (D^*x_0)_g,(D^*w)_g\ra^2}{\|(D^*x_0)_g\|^3}\;\; \mbox{for}\;\;  w\in \dom d^2\Psi(x_0|0). 
\end{eqnarray}

\end{Lemma}
\begin{proof} 
 Since $x_0$ is an optimal solution to \eqref{ITV}, we have $0\in \partial \Psi(x_0)$. With $C=\Phi^{-1}(\Phi x_0)$, we have
\begin{eqnarray*}
d^2\Psi(x_0|0)(w)&=&\liminf_{t\downarrow 0,\; w'\to w}\dfrac{\|D^*(x_0+tw')\|_{\ell_1/\ell_2}+\iota_C(x_0+tw')-\|D^*x_0\|_{\ell_1/\ell_2}-\iota_C(x_0)-t\la 0,w\ra}{\frac{1}{2}t^2}\\
&=&\liminf_{t\downarrow 0,\; w'\to w}\left(\dfrac{\|D^*(x_0+tw')\|_{\ell_1/\ell_2}-\|D^*x_0\|_{\ell_1/\ell_2} }{\frac{1}{2}t^2}+\iota_{\Ker \Phi}(w')\right).
\end{eqnarray*}
It follows that $\dom d^2\Psi(x_0|0)\subset\Ker \Phi$. Since $D^*_{S}x_0=0$, we have
\begin{eqnarray}\label{Q1}\begin{array}{ll}
&\disp\liminf_{t\downarrow 0,\; w'\st{\Ker  \Phi}\to w}\dfrac{\|D^*(x_0+tw')\|_{\ell_1/\ell_2}-\|D^*x_0\|_{\ell_1/\ell_2} }{\frac{1}{2}t^2}\\ 
&\disp=\liminf_{t\downarrow 0,\; w'\st{\Ker  \Phi}\to w}\left(\dfrac{\|D_T^*(x_0+tw')\|_{\ell_1/\ell_2}-\|D_T^*x_0\|_{\ell_1/\ell_2}-t\la De,w'\ra }{\frac{1}{2}t^2}+\dfrac{\la De,w'\ra +\|D_{S}^*w'\|_{\ell_1/\ell_2}}{\frac{1}{2}t}\right).
\end{array}
\end{eqnarray}
Note that the Euclidean norm $\|u\|$ is twice differentiable at $u\neq 0$ with 
\[
\nabla \|u\|=\frac{u}{\|u\|}\qandq \nabla^2\|u\|=\frac{1}{\|u\|}\Id-\frac{1}{\|u\|^3}uu^*
\]
This together with \eqref{smooth} tells us that
 \begin{eqnarray}\label{Q}\begin{array}{ll}
&\disp\liminf_{t\downarrow 0,\; w'\st{\Ker  \Phi}\to w}\dfrac{\|D_T^*(x_0+tw')\|_{\ell_1/\ell_2}-\|D_T^*x_0\|_{\ell_1/\ell_2}-t\la De,w'\ra }{\frac{1}{2}t^2}\\
 &\disp=\disp\liminf_{t\downarrow 0,\; w'\st{\Ker  \Phi}\to w}\sum_{g\in I}\pa{\dfrac{\|(D^*(x_0+tw'))_g\|-\|(D^*x_0)_g\|-t\la e_g,D^*w'\ra }{\frac{1}{2}t^2}}\\
 &=\disp\sum_{g\in I}\dotp{(D^*w)_g}{\dfrac{(D^*w)_g}{\|(D^*x_0)_g\|}-\frac{1}{\|(D^*x_0)_g\|^3}(D^*x_0)_g(D^*x_0)_g^*(D^*w)_g}\\
 &=\disp\sum_{g\in I}\dfrac{\|(D^*w)_g\|^2\|(D^*x_0)_g\|^2-\la (D^*x_0)_g,(D^*w)_g\ra^2}{\|(D^*x_0)_g\|^3}.
\end{array}
 \end{eqnarray}
 
 Since $x_0$ is an optimal solution to \eqref{ITV}, there exists $z\in \R^p$ with $P_Tz=e$, $\|P_{S}z\|_{\ell_\infty/\ell_2}\le 1$ and $Dz\in\Im \Phi^*$. For any $w^\prime \in \Ker \Phi$, we have 
\begin{equation*}
    \la De,w'\ra = \la D_{T}z,w'\ra = -\la D_{S}z,w'\ra=-\la z,D_{S}^*w'\ra.
\end{equation*}
It follows that  $\la De,w'\ra+\|D_{S}^*w'\|_{\ell_1/\ell_2}\ge 0$. Hence we get
\begin{eqnarray*}
\liminf_{t\downarrow 0,\; w'\st{\Ker  \Phi}\to w} \dfrac{\la De,w'\ra+\|D_{S}^*w'\|_{\ell_1/\ell_2} }{\frac{1}{2}t} &=&\begin{cases}
0 & \text{ if } \la De,w\ra+ \|D_{S}^*w\|_{\ell_1/\ell_2}= 0, w\in \Ker \Phi,\\
\infty &\text{ otherwise}.
\end{cases}
\end{eqnarray*}
Since  the ``liminf'' in \eqref{Q} indeed becomes ``lim'', the latter together with \eqref{Q} and  \eqref{Q1} verifies \eqref{dom} and \eqref{d2}. \end{proof}

%and thus combining with \eqref{KerPhi}, \eqref{DT_perp} and by the property of the limit inferior, we have
%\begin{eqnarray}\label{ed2}
 %   d^2\Psi(x_0|0)(w) &\geq& \sum_{r\in I} \left(\dfrac{\|(D^*w)_r\|^2}{\|(D^*x_0)_r\|}-\dfrac{\la (D^*w)_r,(D^*x_0)_r\ra^2}{\|(D^*x_0)_r\|^3}\right)+\iota_{\Ker\Phi}+\iota_{\CC_J(x_0)}.
%\end{eqnarray}
%Thus, $\dom{d^2\Psi(x_0|0)(w)}\subset\CC_J(x_0)\cap \Ker{\Phi}$. Conversely, for any $w\in \CC_J(x_0)\cap \Ker{\Phi}$, by choosing $w'=\lambda w$ with $\lambda=1$, we have $\iota_{\Ker\Phi}(w')=0$ and $\|D_{T^{\perp}}^*w'\|_{\ell_1/\ell_2}+\la De,w'\ra=0$, thus
%%  d^2\Psi(x_0|0)(w)&\leq&\lim_{t\downarrow 0,\; \lambda\to 1}\dfrac{\|D_T^*(x_0+tw')\|_{\ell_1/\ell_2}-\|D_T^*x_0\|_{\ell_1/\ell_2}-t\la De,w'\ra }{\frac{1}{2}t^2}\\
  % &=&\sum_{r\in I} \left(\dfrac{\|(D^*w)_r\|^2}{\|(D^*x_0)_r\|}-\dfrac{\la (D^*w)_r,(D^*x_0)_r\ra^2}{\|(D^*x_0)_r\|^3}\right).
%\end{eqnarray*}
%This together with \eqref{ed2} verifies that   $\dom{d^2\Psi(x_0|0)(w)}=\CC_J(x_0)\cap \Ker{\Phi}$ and \eqref{d2}. The proof is completed. \end{proof} 

This calculation allows us to establish the main result in this section, which gives a quantitative  characterization for unique/strong solutions to $\ell_1/\ell_2$ problem \eqref{ITV}. 

\begin{Theorem}[Characterizations for unique/strong solutions to $\ell_1/\ell_2$ problems]\label{BP12}  The following assertions are equivalent:
\begin{enumerate}[label={\rm (\roman*)}]
    \item $x_0$ is a unique solution to problem \eqref{ITV}.
    \item $x_0$ is a strong solution to problem \eqref{ITV}.
\item $x_0$ is a solution to \eqref{ITV}, $\Ker \Phi\cap \EE\cap \Ker D^*_{S}=\{0\}$,  and 
    \begin{equation}\label{zeta}
    \zeta(e) \eqdef \disp\min_{u\in \Ker  MD_{S}}\|(MD_{S})^\dag MD_Te-u\|_{\ell_\infty/\ell_2}<1,
    \end{equation}
where $M^*$ is a matrix forming a basis matrix to $\Ker \Phi\cap \EE$.

\end{enumerate}
\end{Theorem}
\begin{proof} 
 We first claim that $x_0$ is a solution to \eqref{ITV} if and only if 
\begin{equation}\label{Opt}
    \Ker\Phi \cap \Int ( \CC_{J}(x_0))=\emptyset.
\end{equation}
Indeed, $x_0$ is a solution to \eqref{ITV} if and only if $d\Psi(x_0)(w)\ge 0$ for all $w\in \R^n$. Due to the computation of $d\Psi(\ox)(w)$ in \eqref{dPhi}, $d\Psi(\ox)(w)\ge 0$ means  
\begin{equation}\label{Opt2}
\la De,w\ra+\|D^*_{S}w\|_{\ell_1/\ell_2}\ge 0\quad \mbox{for all}\quad w\in \Ker \Phi,
\end{equation}
which is equivalent to \eqref{Opt}. 

Next let us  verify the equivalence between (i) and (ii). By Theorem~\ref{Uniq}, it suffices to show that condition \eqref{Geo} implies (ii).  Note from \eqref{dom} and \eqref{d2} that if  $d^2\Psi(x_0|0)(w)\le 0$ then  $w\in \Ker \Phi \cap \bd  \CC_{J}(x_0)$ and $w\in \EE.$  Since $\Ker \Phi\cap (\EE\cap \bd  \CC_{J}(x_0))=\{0\}$ by \eqref{Geo}, we have  $d^2\Psi(x_0|0)(w)>0 $ for all $w\neq 0$. It follows from Lemma~\ref{Fa} that $x_0$ is a strong solution. This verifies the equivalence between (i) and (ii).

To justify the equivalence between (i) and (iii), by  \eqref{Opt}, we only need to show that condition 
\begin{equation}\label{CD3}
    \Ker \Phi  \cap \EE\cap \bd  \CC_{J}(x_0)=\{0\}
\end{equation}
 in \eqref{Geo} is equivalent to the combination  of $\Ker\Phi\cap \EE\cap\Ker D^*_{S}=\{0\}$ and \eqref{zeta} provided that  $x_0$ is a solution to \eqref{ITV}. According to \eqref{Opt2},  \eqref{CD3} is equivalent to  the condition that there exists some $c>0$ such that 
\begin{equation}\label{h}
k(w) \eqdef \la De,w\ra+\|D^*_{S}w\|_{\ell_1/\ell_2}\ge c\|w\|\quad \mbox{for all}\quad w\in \Ker\Phi \cap \EE.
\end{equation}
Since $x_0$ is a solution to \eqref{ITV}, there exists $v\in \R^p$ such that $v\in \partial \|D^*x_0\|_{\ell_1/\ell_2}$ and $Dv\in \Im \Phi^*$. It follows that  $D(e+P_{S}v)\in \Ker M=\Im \Phi^*+\EE^\perp$, we have $MD_{S}v=-MDe$. It follows from \eqref{cons} that 
\[
MD_{S}(MD_{S})^\dag MD_Te=MD_Te.
\]
For any $w=M^*u\in 
\Ker \Phi\cap C$ it follows that 
\begin{eqnarray}\label{hM}\begin{array}{ll}
k(w)&=k(M^*u)=\la De,M^*u\ra+\|D^*_{S}M^*u\|_{\ell_1/\ell_2}\\
&=\la MD_Te, u\ra+\|D^*_{S}w\|_{\ell_1/\ell_2}\\
&=\la (MD_{S})^\dag MD_Te,D^*_{S}M^*u\ra+\|D^*_{S}M^*u\|_{\ell_1/\ell_2}.
\end{array}
\end{eqnarray}
Mimicking  the proof of Theorem~\ref{mtheo} by replacing $N$ there by $M$ and $\Ker \Phi$ by its subspace $\Ker \Phi\cap \EE$, the inequality $k(w)\ge c\|w\|$ in \eqref{h} is equivalent to (iii). The proof is complete. \end{proof}

It is worth noting that condition 
\begin{equation}\label{HJC}
    \Ker \Phi\cap \EE\cap \Ker D^*_{S}=\{0\}
\end{equation}
in part (iii) is strictly weaker than  the Restricted Injectivity in Theorem~\ref{mtheo}; see  Example~\ref{Ex12}. It means that $\Phi$ is injective on the subspace $\EE\cap \Ker D^*_{S}$. We refer \eqref{HJC} as {\em Strong Restricted Injectivity} condition. Moreover, we call the constant $\zeta(e)$ {\em Strong Source Coefficient}, while the condition $\zeta(e)<1$ is refered as {\em Analysis Nondegenerate Source Condition} for solution uniqueness to $\ell_1/\ell_2$ problem \eqref{BP0}.

\begin{Remark}[Checking the Strong Restricted Injectivity and the Analysis Nondegenerate Source Condition]
Set the matrix $S \eqdef ((D^*_Tx_0)_g)_{ g\in I}$ to be an $p\times |I|$ matrix, where $|I|$ is the cardinality of $I$.  Observe from  \eqref{d2} that $w\in \EE$ if and only if the following system 
 \[
 D^*_Tw=S\lm \mbox{ is consistent  with }\lm\in \R^{|I|}.
 \]
 Since $S$ is injective, we have $\lm=S^\dag  D^*_Tw =(S^*S)^{-1}S^*D^*_Tw$ due to \eqref{cons}. It follows that $w\in \EE$ if and only if \begin{equation}\label{Si}
 S(S^*S)^{-1}S^*D^*_Tw-D^*_Tw=0
 \end{equation}
Note that  $S^*S={\rm diag}\, (\|(D^*_Tx_0)_g\|^2)_{b\in I}$ is an $|I|\times |I|$ diagonal matrix. Representing $S$ in terms of groups, we have 
 \[
 S(S^*S)^{-1}S^*={\rm diag\,}\left(\sum_{g\in I}\|(D_T^*x_0)_g\|^{-2}\delta_g\right)SS^*,
 \]
 where $\delta_g=(0,\ldots,0, \underbrace{1,1,\ldots,1}_\text{in $V_g$},0,\ldots,0)$ is the unit vector in $V_g$. This together with \eqref{Si} tells us that
 $\EE$ is the kernel of the following matrix
\begin{equation}\label{Q3}
     Q=\diag\left(\sum_{g\in I}\|(D_T^*x_0)_g\|^2\delta_g\right)D_T^*-SS^*D_T^*\quad\mbox{with}\quad S=((D^*_Tx_0)_g)_{ g\in I}.
 \end{equation}
The Strong Restricted Injectivity \eqref{HJC} is equivalent to $\Ker\begin{pmatrix}\Phi\\Q\end{pmatrix}\cap \Ker D^*_{S}=\{0\}. $

Furthermore, $M^*$ forms a basis matrix of $\Ker \begin{pmatrix}\Phi\\Q\end{pmatrix}$, which is found from the SVD of $\begin{pmatrix}\Phi\\Q\end{pmatrix}$. Similarly to \eqref{rhoz}, $\zeta(e)$ is the optimal solution to 
\begin{equation}\label{Comz}
\min \|z\|_{\ell_2/\ell_\infty}\qstq MDz=-MDe\qandq z\in \bigoplus_{g\in I}V_g. 
\end{equation}
So $\zeta^2(e)$ is the optimal value to the following convex optimization problem 
\begin{equation}\label{zeta2}
\min_{t \geq 0,z}  t \qstq  MDz=-MDe, \quad \|z_g\|^2-t\le 0,\quad g\in I\quad  \mbox{and}\quad z\in\bigoplus_{g\in I}V_g 
\end{equation}
with $|I|G+1$ variables,  which can be solved by available packages such as \texttt{cvxopt}; see Section 6 for further discussion.
\end{Remark}

Next we show that $\zeta(e)\le \rho(e)$ when $x_0$ is an optimal solution to \eqref{ITV}. 

\begin{Proposition}[Comparison between $\rho(e)$ and $\zeta(e)$]\label{Com} Suppose that $x_0$ is an optimal solution to \eqref{ITV}. Then we have $\zeta(e)\le \rho(e)\le 1$.
\end{Proposition}
\begin{proof} 

It follows from \eqref{yz} that 
\[
ND_{S}(ND_{S})^\dag ND_Te=ND_Te\qandq \rho(e)\le 1. 
\]
For any $w\in \Ker\Phi \cap\EE=\Im M^*\subset \Ker\Phi=\Im N^*$, it is similar to \eqref{hM} that
\[
\la De,w\ra=\la (MD_{S})^\dag MD_Te,D^*_{S}w\ra=\la (ND_{S})^\dag ND_Te,D^*_{S}w\ra
\]
Note from the definition in \eqref{zeta} that 
\begin{eqnarray}\label{zen}\begin{array}{ll}
\zeta(e)&=\disp\sup_{\|v\|\le 1}\la -(MD_{S})^\dag MD_Te,v\ra-\iota_{\Im D^*_{S}M^*}(v)\\
&= \disp\sup_{\|v\|\le 1}\la -(ND_{S})^\dag ND_Te,v\ra-\iota_{\Im D^*_{S}M^*}(v)\\
&=\disp\min_{u\in \Ker  MD_{S}}\|(ND_{S})^\dag ND_Te-u\|_{\ell_\infty/\ell_2}.
\end{array}
\end{eqnarray}
Since $\Ker ND_{S}\subset \Ker  MD_{S}$, we have $\zeta(e)\le \rho(e)$. \end{proof}

Although $\zeta(e)$ could be computed by involving $N$ via \eqref{zen}, the format in \eqref{zeta} is more preferable. This is due to the fact that the Moore-Penrose inverse $MD_{S}^\dag$ has a closed form as $(MD_{S})^* (MD_{S}(MD_{S})^*)^{-1}$  when the Strong Restricted Injectivity  \eqref{HJC} is in charged. In general, $\rho(e)$ is strictly smaller $\zeta(e)$. This fact is obtained through numerical experiments in our Section~6.

 By replacing $\zeta(e)$ by $\rho(e)$ in Theorem~\ref{BP12}, we do not need to assume $x_0$ to be an optimal solution, but we only have a sufficient condition for solution of uniqueness. %Note that this condition is not necessary as discussed in our Example~\ref{Ex12}.  
\begin{Corollary}[Sufficient condition for  solution uniqueness to $\ell_1/\ell_2$ problem]\label{CoBP} $x_0$ is the unique solution to problem \eqref{ITV} provided that $\rho(e)<1$ and both conditions $\eqref{yz}$ and \eqref{HJC} are satisfied. 
\end{Corollary}
\begin{proof} 
 It follows directly from Proposition~\ref{Sol}, Theorem~\ref{BP12}, and Proposition~\ref{Com}. \end{proof}

 An simple upper bound for $\zeta(e)$ is
\begin{equation}\label{gamma}
\zeta(e)\le \gg(e) \eqdef \|(MD_{S})^\dag MD_Te\|_{\ell_\infty/\ell_2},
\end{equation}
which is also used in Section 6 to check solution uniqueness. The inequality is indeed strict. The following result  is straightforward from Theorem~\ref{BP12}.

\begin{Corollary}[Sufficient condition for solution uniqueness to $\ell_1/\ell_2$ problem]\label{Corogg} Suppose that $x_0$  is an optimal solution to \eqref{ITV}. Then the Strong Restricted Injectivity \eqref{HJC}  at $x_0$ and $\gg(e)<1$ 
are sufficient for solution uniqueness at $x_0$. They become necessary conditions provided that $MD_{S}$ is injective.  
\end{Corollary}
%\begin{proof} 
 %The result follows from Theorem~\ref{BP12} by noting that 
%\[
%\zeta(e)\le \|(MD_{S})^\dag MD_Te\|_{\ell_\infty/\ell_2}.
%\]
%The equality occurs when if $MD_{S}$ is injective. \end{proof} 

In the spirit of Proposition~\ref{Com}, it is possible that $x_0$ is a unique solution to \eqref{ITV} but $\rho(e)=1$. It means that the Nondegenerate Source Condition may not happen. However, the following result provides characterizations for solution uniqueness whose statement closely relates to the Nondegenerate Source Condition.

\begin{Corollary}[Characterization for solution uniqueness to problem \eqref{ITV}]\label{LoBP} The following are equivalent:
\begin{enumerate}[label={\rm (\roman*)}]
    \item $x_0$ is the unique solution to problem \eqref{ITV}.

\item With an arbitrary $v\in \R^p$ satisfying $Dv\in \Im \Phi^*$, $P_Tv=e$, and $\|P_{S}v\|_{\ell_\infty/\ell_2}\le 1$,  the following system has only trivial solution
\[
\Phi w=0, \quad w\in \EE, \qandq D_{S}^*w\in \cone \{v_g|g\in L\}\times\{0_H\},
\]
where $L \eqdef \{g\in K|\; \|v_g\|=1\}$ and $H \eqdef K\setminus L.$

\end{enumerate}
\end{Corollary}
\begin{proof} 
 Pick an arbitrary $v\in \R^p$ satisfying $Dv\in \Im \Phi^*$, $P_Tv=e$, and $\|P_{S}v\|_{\ell_\infty/\ell_2}\le 1$. Such an $v$ always exists as when $x_0$ is a solution to \eqref{ITV}, i.e., $\Ker \Phi \cap \Int \CC_{J}(x_0)=\emptyset$. We claim that 
\begin{equation}\label{D2}
    \Ker \Phi\cap \bd \CC_{J}(x_0)=\{w\in \Ker\Phi |\; D_{S}^*w\in \cone \{v_g|g\in L\}\times\{0_H\}\} 
    \end{equation}
with   $L=\{g\in K|\; \|v_g\|=1\}$ and $H=K\setminus L.$ For any $w\in \Ker\Phi$, note that 
\[
\begin{array}{ll}
\disp\la De,w\ra+\|D^*_{S}w\|&=\la Dv,w\ra-\la DP_{S}v,w\ra+\|D^*_{S}w\|\\
&=\disp -\la P_{S}v,D^*_{S}w\ra+\| D^*_{S}w\|_{\ell_1/\ell_2}\\
&\disp=\sum_{g\in L}\left(-\la v_g, (D^*_Tw)_g\ra+\|(D^*_Tw)_g\|\right)+\sum_{g\in H}\left(\la v_g, (D^*_Tw)_g\ra+\|(D^*_Tw)_g\|\right)\\
&\disp=\sum_{g\in L}\left(-\la v_g, (D^*_Tw)_g\ra+\|(D^*_Tw)_g\|\right)+\sum_{g\in H}\|(D^*_Tw)_g\|\ge 0.
\end{array}
\]
It follows that $w\in\bd \CC_{J}(x_0)$ if and only if there exist $\lm_g\ge 0 $, $g\in L$ such that $(D^*_Tw)_g=\lm_gv_g$ for any $g\in L$ and $(D^*_Tw)_g=0$ for any $g\in H$, which verifies \eqref{D2}. The equivalence between (i) and (ii) follows from   Theorem~\ref{Uniq}. \end{proof}

%Although Corollary~\ref{LoBP} provides another characterization for solution uniqueness, using it in practice is not simple. After solving \eqref{ITV} we have 

As $\Ker\Phi\cap\EE\cap \bd \CC_{J}(x_0)\subset\Ker\Phi\cap\Ker D^*_{V_H}$ with  $\disp V_H \eqdef \bigoplus_{g\in H}V_g$, the following result is straightforward from Corollary~\ref{CoBP}.
\begin{Corollary}[Sufficient condition to the solution uniqueness to problem \eqref{ITV}]\label{SufCoBP} $x_0$ is the unique solution to problem \eqref{ITV} provided that there exists some $v$ such that $Dv\in \Im \Phi^*$, $P_Tv=e$,  $\|P_{S}v\|_{\ell_\infty/\ell_2}\le1$, and that  
\begin{equation}\label{H}
\Ker\Phi\cap\Ker D^*_{V_H} =\{0\} 
\end{equation}
with $K=\{g\in J|\; \|v_g\|=1\}$,  $H=J\setminus K$, and $\disp V_H \eqdef \bigoplus_{g\in H}V_g$.
\end{Corollary}

 The sufficient condition in \cite[Proposition~7.1]{G11} is a special case of Corollary~\ref{SufCoBP} where $D^*$ is the discrete gradient operator. Moreover, \cite[Theorem~3.4]{JKL15} and \cite[Theorem~3]{RF08} even assume a stronger condition as they require $K=\emptyset$. Let us revisit Example~\ref{Ex12}. Pick any $v\in \partial \|x_0\|_{\ell_1,\ell_2}\cap \Im\Phi^*$, i.e., $v\in (0,1)\times [-1,1]\cap \Im\{(1,1,0),(1,0,-1)\}$. It follows that $v=(0,1,1)$ and thus $K=\{3\}$, $H=\emptyset$, which clearly implies \eqref{H}.

%In all the results about characterizations for solution uniqueness to group-sparsity problem \eqref{ITV} up till now, we need to assume in advance that $x_0$ is an optimal solution in Theorem~\ref{BP12}. These results are particularly helpful when after solving \eqref{ITV} we have a solution and then using them to check the uniqueness. However, it is possible to characterize solution uniqueness to group-sparsity problem \eqref{BP0} without assuming initially about optimality as follows.

According to Theorem~\ref{BP12}, Theorem~\ref{SRob}, and Proposition~\ref{NRC}, solution uniqueness to group-sparsity problem \eqref{ITV} is equivalent to the robust recovery with rate $\mathcal{O}(\sqrt{\delta})$.

\begin{Corollary}[Robust recovery and solution uniqueness  for group-sparsity problems] The following statements are equivalent: 
\begin{enumerate}[label={\rm (\roman*)}]
\item $x_0$ is a unique solution to \eqref{ITV}.

\item For sufficiently small $\delta>0$, any solution $x_\delta$ to problem \eqref{BP1} with $J_0(\cdot)=\|\cdot\|_{\ell_1/\ell_2}$ satisfies $\|x_\delta-x_0\|\le \mathcal{O}(\sqrt{\delta})$ whenever $\|y-y_0\|\le \delta.$

 \item For any $c_1>0$ and sufficiently small $\delta>0$, any solution $x_\mu$ to \eqref{Las3} with $J_0(\cdot)=\|\cdot\|_{\ell_1/\ell_2}$ satisfies $\|x_\mu-x_0\|\le \mathcal{O}(\sqrt{\delta})$ whenever $\|y-y_0\|\le \delta$ and $\mu=c_1\delta$.
\end{enumerate}
\end{Corollary}
%\textcolor{red}{
%\begin{Remark}{\rm 
%\cite{V14} showed that Non-degenerate Source Condition (NDSC) and Restrictive Injectivity Condition (INJ) together imply the solution uniqueness of $x_0$ and the linear convergence rate of $x_{\mu}$ to $x_0$. It is worth to emphasize that our result in Corollary~\ref{LoBP} does not contradict \cite{V14}. In fact, the conditions (NDSC) and (INJ) are too strong for the solution uniqueness, which is showed in Example~\ref{Ex12} and Example~\ref{Ex52}.   }
% \end{Remark}
%}

From Theorem~\ref{BP12}, it is natural to raise the following question for other decomposable norm minimization problem. 
\begin{enumerate}[label={\rm \bf (Q.\arabic*)},start=3]
\item If $x_0$ is a unique solution to \eqref{BP0}, can it be a strong solution to \eqref{BP0}? \label{OQ3}
\end{enumerate}
The answer to \ref{OQ3} is affirmative for $\ell_1/\ell_2$ problem as in Theorem~\ref{BP12}. For the $\ell_1$ problem when $\|\cdot\|_\AA=\|\cdot\|_1$, we also have the positive answer due to Corollary~\ref{ell2}. However, it is not the case for the  nuclear norm minimization problem. The following example modifies \cite[Example~3.1]{BLN19} to prove that claim.

\begin{Example}[Difference between unique solution and strong solution to NNM]\label{Last}  
Consider the following optimization problem 
\begin{equation}\label{Impo}
\min_{X\in\R^{2\times 2}} \|X\|_* \qstq \Phi(X) \eqdef 
\begin{pmatrix}
     X_{11}+X_{22}  \\
     X_{12}-X_{21}+X_{22} 
\end{pmatrix}=
\begin{pmatrix}
     1  \\
     0
\end{pmatrix}.
\end{equation}
For any $X=\begin{pmatrix}
     a &b  \\
     c &d
\end{pmatrix}$ with $a+d=1$ and $b-c+d=0$, we obtain from \eqref{nucl} that 
\[
\|X\|_*=\sqrt{a^2+b^2+c^2+d^2+2|ad-bc|}\ge \sqrt{|a+d|^2+|b-c|^2}\ge 1,
\]
where the equality occurs when $\|X\|_*=1$, $b-c=0$, $a+d=1$, and $b-c+d=0$,  which means $b=c$, $a=1$, $d=0$, and  $\|X\|_*=\sqrt{1+b^2+c^2+2|bc|}=1$. So $\Bar X=\begin{pmatrix}
     1 &0  \\
     0 &0
\end{pmatrix}$ is the unique solution to problem \eqref{Impo}. Choose $X_\ve=\begin{pmatrix}
     1-\ve^{1.5}&\ve-\ve^{1.5}   \\
     \ve &\ve^{1.5}
\end{pmatrix}$ with $\ve>0$ sufficiently small and note that $X_\ve$ satisfies the equation in \eqref{Impo}. It follows that 
\begin{align*}
\|X_\ve\|_*-\|\Bar X\|_*&=\disp\sqrt{(1-\ve^{1.5})^2+(\ve-\ve^{1.5})^2+\ve^2+\ve^3+2|(1-\ve^{1.5})\ve^{1.5}-(\ve-\ve^{1.5})\ve|}-1\\
&=\disp\sqrt{(1-\ve^{1.5})^2+(\ve-\ve^{1.5})^2+\ve^2+\ve^3+2(1-\ve^{1.5})\ve^{1.5}-2(\ve-\ve^{1.5})\ve}-1\\
&=\disp\sqrt{1+\ve^3}-1={\mathcal O}(\ve^3) .
\end{align*}
Moreover, $\|X_\ve-\Bar X\|_F^2=\ve^3+(\ve-\ve^{1.5})^2+\ve^2+\ve^3=O(\ve^2)$. This tells us that $\Bar X$ is not a strong solution to \eqref{Impo}. 
\end{Example}

%%%%%%%%%%%%%%%%%%%%%%%%%%%
\subsection{Connections between unique/strong solutions in the noiseless case}
Let us consider the particular case of problem \eqref{Las0}, which reads
\begin{equation}\label{Las12}
    \min \Theta(x) \eqdef f(\Phi x)+\mu\|D^*x\|_{\ell_1/\ell_2}
\end{equation}
with constant $\mu>0$. In the following result, we show that  a unique solution to \eqref{Las12} is also a strong solution.  Consequently, all the characterizations for solution uniqueness to \eqref{ITV} in this section can be used to characterize unique/strong solutions to problem \eqref{Las12} due to Proposition~\ref{Stheo}. Moreover, this result gives an affirmative answer for \ref{OQ2} in the previous section.

\begin{Theorem}[Characterization to solution uniqueness to $\ell_1/\ell_2$ regularized problem]\label{Las2} Suppose that $\ox$ is an optimal solution to problem \eqref{Las12}. The following are equivalent:
\begin{enumerate}[label={\rm (\roman*)}]
\item $\ox$ is the unique solution to \eqref{Las12}.

\item $\ox$ {is} the strong solution to \eqref{Las12}.

\item $\ox$ is the strong solution to \eqref{ITV} with $x_0=\ox$.
\end{enumerate}

\end{Theorem}
\begin{proof} 
 (i) and (iii) are equivalent due to Proposition~\ref{Stheo} and Theorem~\ref{BP12}. [(ii)$\Rightarrow$(i)] is trivial. It suffices to verify [(iii)$\Rightarrow$(ii)]. Suppose that $\ox$ is a strong solution to \eqref{ITV}. By Theorem~\ref{BP12} and Theorem~\ref{Uniq}, $\Ker \Phi\cap \EE\cap \bd \CC_{J}(x_0)=\{0\}$.  Since $\ox$ is an optimal solution to \eqref{Las2}, we have $$0\in \partial \Theta(\ox)=\Phi^*\nabla f(\Phi \ox)+D\partial\|D^*\ox\|_{\ell_1/\ell_2}\subset D\partial\|D^*\ox\|_{\ell_1/\ell_2}+\Im\Phi^*.$$ Thus there exists $\oz\in \partial\|D^*\ox\|_{\ell_1/\ell_2}$ such that $D\oz=-\Phi^*\nabla f(\Phi \ox)$, $P_T\oz=e$, and $\|P_{S}\oz\|_{\ell_\infty/\ell_2}\le 1$, where $e,T$, and $I$ are defined at the beginning at the section with $x_0=\ox$.

For any $w\in \R^n$, note further that 
\begin{equation}\label{Stro}
d^2\Theta(\ox|0)(w)=\la \Phi w,\nabla^2f(\Phi\ox)\Phi w\ra+d^2 J(\ox|\; D\oz)(w)
\end{equation}
with $J(x)=\|D^*x\|_{\ell_1/\ell_2}$. Since $D^*_{S}\ox=0$, we have
\[\begin{array}{ll}
&d^2 J(\ox|\; D\oz)(w)\disp=
\liminf_{t\downarrow 0,\; w'\to w}\dfrac{\|D^*(\ox+tw')\|_{\ell_1/\ell_2}-\|D^*\ox\|_{\ell_1/\ell_2}(\ox)-t\la D\oz,w'\ra}{\frac{1}{2}t^2}\\
&\disp=\liminf_{t\downarrow 0,\; w'\to w}\left(\dfrac{\|D_T^*(x_0+tw')\|_{\ell_1/\ell_2}-\|D_T^*x_0\|_{\ell_1/\ell_2}- t\la De,w'\ra }{{\frac{1}{2}t^2}}+\dfrac{{-}\la D_{S}\oz,w'\ra {+} \|D_{S}^*w'\|_{\ell_1/\ell_2}}{\frac{1}{2}t}\right).
\end{array}
\]
Note further that 
\[
-\la D_{S}\oz,w'\ra +\|D_{S}^*w'\|_{\ell_1/\ell_2}= -\la P_{S}\oz,D_{S}^*w'\ra+\|D_{S}^*w'\|_{\ell_1/\ell_2}\ge 0
\]
as $\|P_{S}\oz\|_{\ell_\infty/\ell_2}\le 1$. 
It is similar to \eqref{Q1} and \eqref{Q} that 
\begin{eqnarray}
\dom d^2 J(\ox|\; D\oz)&=&\{w\in \R^n|\; - \la DP_{S}\oz,w\ra +\|D_{S}^*w\|_{\ell_1/\ell_2}=0\} \eqdef \mathcal{F},\label{d21}\\
d^2 J(\ox|\; D\oz)(w)&=&\sum_{g\in I}\dfrac{\|(D_T^*w)_g\|^2}{\|(D_T^*x_0)_g\|}-\dfrac{\la (D_T^*w)_g,(D_T^*x_0)_g\ra^2}{\|(D_T^*x_0)_g\|^3}.\label{d22}
\end{eqnarray}
Since $\nabla^2f(\Phi \ox)\succ 0$, we derive from \eqref{Stro}, \eqref{d21}, and \eqref{d22} that  
\begin{equation*}
\Ker d^2\Theta(\ox|0)=\Ker \Phi\cap \EE\cap \mathcal{F}. 
\end{equation*}
Furthermore, with  $w\in \Ker \Phi$ we have $\Phi w=0$ and 
\[
\la DP_{S}\oz,w\ra=\la D\oz,w\ra-\la De,w\ra=\la -\Phi^*\nabla f(\Phi \ox),w\ra-\la De,w\ra=-\la De,w\ra.
\]
This together with \eqref{d21} implies that  $$\Ker d^2\Theta(\ox|0)=\Ker \Phi\cap \EE\cap \mathcal{F}=\Ker \Phi\cap \EE\cap \bd \CC_{J}(x_0)=\{0\}.$$ By Lemma~\ref{Fa}, $\ox$ is a strong solution to problem \eqref{Las12}. The proof is complete.  \end{proof}

%% file: tex/sec_numerics.tex
%%%%%%%%%%%%%%%%%%%%%%%%%%%%%%%%%%%%%%%%%%%%%%%%%%%%%%%%
\section{Numerical verification of solution uniqueness for group-sparsity}
\label{sec:numerics}
%%%%%%%%%%%%%%%%%%%%%%%%%%%%%%%%%%%%%%%%%%%%%%%%%%%%%%%%
With an aim of demonstrating that our conditions for sharp, strong/unique solution are verifiable, we have implemented a simulation using synthetic data. In our simulation study, $\Phi$ was generated as an $(2000\times 240)$ Gaussian matrix whose entries are independently and identically drawn from the standard normal distribution, $\NN(0,1)$. We next randomly divided the set of indicators range from $1$ to $2000$ into $k=100$ groups of size $20$ with $5$ randomly selected active groups. Then, we constructed a measured signal of length $m=240$, $y_0 \eqdef \Phi x_0$, based on the original signal $x_0$ whose elements in each active group are independently and identically distributed $\NN(0,1)$. We then used our proposed conditions to verify whether $x_0$ is a solution to \eqref{BP0} using the conditions in Proposition~\ref{Sol} and summarize the number of cases where $x_0$ is classified as sharp or unique/strong solution by the criteria from Theorem~\ref{mtheo} and \ref{BP12}.
For checking strong and sharp minima, we only need to compute the Source Coefficient $\rho(e)$ and the Strong Source Coefficient $\zeta(e)$ whenever $\tau(e)$ in \eqref{xi} or $\gamma(e)$ in \eqref{gamma} is greater than or equal to $1$, respectively, since calculating these numbers are much easier. Similar to the scheme for calculating $\zeta(e)$ in \eqref{zeta2}, $\rho^2(e)$ is the optimal value to the following convex optimization problem (recall \eqref{rhoz})
\begin{equation}\label{rho2}
\min_{t \geq 0, z} t \qstq ND z=-NDe,\quad \|z_g\|^2\le t, g\in K,\qandq z\in \bigoplus_{g\in K}V_g. 
\end{equation}
where $K$ is the set of inactive groups.  Note that \eqref{zeta2} and \eqref{rho2} are second-order cone programming problems and can be solved via function \texttt{solvers.socp} of \texttt{cvxopt} package. In our experiment, $\rho(e)$ or $\zeta(e)$ are calculated when $\tau(e)$ or $\gamma(e)$ is greater than or equal to $0.99$.

The results are recorded in the following tree diagram.
\begin{center}
\begin{tikzpicture}[sibling distance=5cm,
  every node/.style = {shape=rectangle, rounded corners,
    draw, align=center,
    top color=white, bottom color=white!20}]]
 \node {$x_0$ is a solution \\$100$ tests}
    child { node {$\Ker\Phi\cap\Ker D^*_{S} =\{0\}$\\ 100 tests} 
        child { node {$\tau(e)<0.99$\\ 11 cases}
           }
        child { node {$\tau(e)\geq 0.99$\\ 89 cases} 
            child { node {$\rho(e)<0.95$\\ 60 cases} }
            child { node {$0.95\leq\rho(e)<1.05$\\29 cases} 
                child { node {$\gamma(e)<0.99$\\ 26 cases} } 
                child { node {$\gamma(e)\geq 0.99$\\ 3 cases} 
                    child { node {$\zeta(e)<0.95$\\ 3 cases} }}}}
    };
\end{tikzpicture}
\end{center}
In all $100$ tested random cases, $x_0$ is verified as a solution to \eqref{BP0} and satisfying the Restricted Injectivity Condition. Among them, there are $11$ cases with $\tau(e)<0.99$ thus are classified as sharp solution, the rest are passed to next step for calculating $\rho(e)$. There are $60$ cases with $\rho<0.95$ and $29$ tests with $0.95<\rho(e)<1.05$. Hence, we had $71$ cases $x_0$ is the sharp solution. We continue the experiment by checking the strong solution condition on the rest $29$ cases. Note that since all cases satisfy the Restricted Injectivity Condition, they automatically satisfy the Strong Restricted Injectivity Condition and it is left to check the Analysis Nondegenerate Source Condition. All $29$ cases are indicated as satisfying the Analysis Nondegenerate Source Condition with $26$ cases having $\gamma(e)<0.99$ and $3$ cases where $\gamma(e)\geq 0.99$ and $\zeta(e)<0.95$. It means that we have 29 cases of strong solutions that are not sharp. 
\begin{table}[h!]
    \centering
    \begin{tabular}{l c}
    \hline
   & number of cases \\
     \hline

\medspace
Sharp solution & 71 \\

\medspace
Strong solution (non-sharp) & 29\\
\hline

    \end{tabular}
    \caption{Number of cases with strong and sharp solutions}
    \label{tab:strong}
\end{table}

% \begin{table}[h!]
%     \centering
%     \begin{tabular}{|c|c|c|c|}
%     \hline
%      & scenario 1 & scenario 2 & scenario 3 \\
%      \hline
% unique solution & 87 & 100 & 100\\
% \hline
% $\Ker\Phi\cap\EE\cap\Ker D_{T^{\perp}} =\{0\}$ & 87 & 100 & 100\\ 
% \hline
% $\gamma(e)<0.9$ & 81 & 99 & 100\\
% \hline
% $0.9\leq \gamma(e)<1$ & 6 & 1 & 0 \\
% \hline
% $\gamma(e)\geq 1$ & 0 & 0 & 0\\
% \hline
%     \end{tabular}
%     \caption{Number of cases that have unique solutions}
%     \label{tab:unique}
% \end{table}

%% file: tex/sec_conclusion.tex
%%%%%%%%%%%%%%%%%%%%%%%%%%%%%%%%%%%%%%%%%%%%%%%%%%%%%%%%
\section{Conclusion}
\label{sec:conclusion}
%%%%%%%%%%%%%%%%%%%%%%%%%%%%%%%%%%%%%%%%%%%%%%%%%%%%%%%%
In this paper we show that sharp minima and strong minima play important roles for robust recovery with different rates. We also provide some quantitative  characterizations for sharp solutions to convex regularized problems. Unique solutions to $\ell_1$ problems are actually sharp solutions. For group sparsity problems, unique solutions are strong solutions. We also obtain several conditions guaranteeing solution uniqueness to group sparsity problems.   As solution uniqueness to $\ell_1$ problem plays a central role in the area of exact recovery with high probability, we plan to use our results to find a better bound for exact recovery for group-sparsity problems in comparison with the one obtained in \cite{CR13,RRN12}, at which they only use sufficient conditions for solution uniqueness. 

Example~\ref{Last3} and Example~\ref{Last} raise important questions about the solution uniqueness and strong minima for the nuclear norm minimization problems. Unique solutions to \eqref{BP0} with the nulear norm are neither sharp nor strong solutions. But a strong solution in this case is certainly a unique solution. It means that second-order analysis can provide a sufficient condition for solution uniqueness, and such a condition should be weaker than the one in Theorem~\ref{mtheo} for sharp solutions. However, unlike the analysis in Lemma~\ref{SOD}, the second subderivative of the nuclear norm is far more intricate to compute. Understanding  solution uniqueness and strong minima for the case of the nuclear norm, or more generally for spectral functions, is a project that we plan to pursue in the future. %We need geometric characterization for strong solution like (iii) in Theorem~\ref{BP12} to proceed further.  